\documentclass[11pt]{amsart}
\usepackage{amsthm}
\usepackage{amsmath}
\usepackage{amssymb}
\usepackage{latexsym}
\usepackage{bm}
\usepackage{pstricks}
\usepackage{pst-plot}
\usepackage{psfrag}
\usepackage{color, graphicx}
\usepackage{ setspace}
\usepackage[all]{xy}
\usepackage[margin=1.25in]{geometry}
\theoremstyle{plain}
\newtheorem{thrm}{Theorem}[section]
\newtheorem{prop}[thrm]{Proposition}
\newtheorem{lem}[thrm]{Lemma}
\newtheorem{cor}[thrm]{Corollary}

\theoremstyle{definition}
\newtheorem{defn}[thrm]{Definition}
\theoremstyle{remark}
\newtheorem{remark}[thrm]{Remark}
\newtheorem{example}[thrm]{Example}

\theoremstyle{plain}
\newtheorem*{intro isom lcc}{Proposition~\ref{isom lcc}}
\newtheorem*{intro ext. induces isom.}{Proposition~\ref{ext. induces isom.}}
\newtheorem*{intro n-conc iff nth invt}{Theorem~\ref{n-conc iff nth invt}}
\newtheorem*{intro ka is conc invt}{Corollary~\ref{ka is conc invt}}
\newtheorem*{intro dense J}{Theorem~\ref{dense J}}
\newtheorem*{intro char has triv}{Theorem~\ref{char has triv}}

\newcommand\Z{\mathbb{Z}}
\newcommand\R{\mathbb{R}}

\newcommand\M{G}

\newcommand\PK{\pi_K}

\newcommand\PJ{\pi_J}

\newcommand\PV{\pi_V}
\newcommand\PC{\pi_C}
\newcommand\HK{\G \g_K}
\newcommand\HJ{\G \g_J}
\newcommand\HV{\G \g_V}
\newcommand\HVn{{\Gamma_n \gamma_V}}
\newcommand\HKn{\G_n \g_K}
\newcommand\HJn{\G_n \g_J}
\newcommand\HJnn{\G_{n+1} \g_J}

\newcommand\HJj{\G_{j} \g_J}
\newcommand\HKj{\G_{j} \g_K}
\newcommand\HVj{\G_{j} \g_V}
\newcommand\HJjj{\G_{j+1} \g_J}
\newcommand\HKjj{\G_{j+1} \g_K}
\newcommand\HVjj{\G_{j+1} \g_V}

\newcommand\EK{{E_K}}
\newcommand\EJ{{E_J}}
\newcommand\EV{{E_V}}

\newcommand\mK{\mu_K}
\newcommand\mJ{\mu_J}

\newcommand\ka{\mathfrak{h}}

\newcommand\g{\gamma}
\newcommand\gK{\gamma_K}
\newcommand\gJ{\gamma_J}

\newcommand\gV{\gamma_V}
\newcommand\G{\Gamma}

\title{Homotopy properties of knots in prime manifolds}
\author{Prudence Heck}
\date{\today}                                           

\begin{document}

\begin{abstract}
We define homotopy-theoretic invariants, $\ka$-invariants, of knots in prime manifolds.  Fix a knot $J$ in a prime manifold $M$.  Call a knot $K \subset M$ concordant to $J$ if it cobounds a properly embedded annulus with $J$ in $M\times I$, and call $K$ $J$-characteristic if there is a degree-one map $\alpha: M \to M$ throwing $K$ onto $J$ and mapping $M-K$ to $M-J$.  These $\ka$-invariants are invariants of concordance and of $J$-characteristicness when $\alpha$ induces the identity on $\pi_1(M)$, and may be viewed as extensions of Milnor's $\overline{\mu}$-invariants.  We do not require the knots considered here to be rationally null-homologous or framed.
\end{abstract}

\maketitle

\section{Introduction}


It is well known that closed, orientable,  prime 3-manifolds come in three forms:  the first type have finite fundamental group and universal cover $S^3$, the second type have infinite, noncyclic fundamental group and contractible universal cover, and the third type is $S^1 \times S^2$.  Manifolds of the first two types are irreducible, meaning that every embedded 2-sphere bounds a ball.  We define invariants for any knot $K$ in a closed, orientable, prime 3-manifold $M$; no choice of framing for $K$ is necessary.  However, we orient $K$ when it represents a nontrivial class in $\pi_1(M)$.\\  



We consider concordance of knots in $M$, where knots $K$ and $J$ are {\em concordant} in $M$ if they cobound a properly embedded cylinder in $M \times I$.  Concordant knots are freely homotopic, making the initial obstruction to concordance of knots the difference of their free homotopy classes.  Yet deeper homotopy theoretic obstructions underlie concordance, and we investigate these here.  We define new homotopy theoretic concordance invariants for knots in $M$.  These invariants extend to knots in three manifolds the prior obstructions defined in ~\cite{milnor1, massey, orr2, orr1}.  The invariants defined here are invariants of the {\em Poincar\'e embedding type} of the knot.  The term, while known to a few specialists, to my knowledge has not appeared in the literature in this context, so we elaborate below.\\

Let $\gamma \in \pi_{1}(M)$ be a fixed homotopy class.  A knot $K$ in the conjugacy class of $\gamma$ determines an inclusion of it's sphere bundle in $M$, $\alpha_K: T^{2} \hookrightarrow \EK$.  If $K$ and $J$ are concordant via a concordance $C$ and if $E_{C} \subset M\times I$ is the exterior of $C$ then the images of these homotopy classes, $[\alpha_K]$ and $[\alpha_J]$, agree in the set $[T^{2}, E_{C}]$.  Although the complement of a knot is not a concordance invariant, the knot determines the homology type of its exterior.  In particular, as observed in Cappell and Shaneson~\cite{cappell-shaneson-codim2}, the exterior of a concordance is a $\Z[\pi_{1}(M)]$-homology cobordism between $E_K$ and $E_J$.  Vogel defined a notion of homology localization of spaces, a functor $\wedge$ on pointed CW-complexes which is initial among functors that turn continuous maps between finite complexes with contractible cofiber into homotopy equivalences, together with a natural transformation from the identity functor to $\wedge$, ~\cite{vogel, ledimet}.  In particular, a homology equivalence $f: Y \hookrightarrow X$ inducing a normally surjective homomorphism on fundamental groups induces a homotopy equivalence $\hat{f}: \widehat{Y} \to \widehat{X}$.  Returning to concordance, it follows that the image of the homotopy class $[\alpha_K]$ under the set map $[T^2, \EK] \to [T^{2}, \widehat{E_K}]$ is a homotopy theoretic concordance invariant of $K$.  This invariant, a homological analogue for codimension two embeddings of the notion of Poincar\'e embedding type that appears in Wall ~\cite{Wall}, first appeared in the context of concordance of disk links in work of Le Dimet ~\cite{ledimet}.  Our invariants are invariants of the Poincar\'e embedding type of the knot in this sense.  However, localization of spaces is an obscure construction, for which very little is known and for which computations are currently out of reach.  By contrast, we expect our invariants to be highly computable, reducing to computations in three dimensional group homology. \\

After fixing a Poincar\'e embedding type one can try to construct a concordance between two knots by choosing a cobordism rel boundary of the knot exteriors, and using surgery techniques to modify this cobordism to a homology cobordism.  This idea underpins Cappell and Shaneson's approach to high-dimensional codimension-two embeddings ~\cite{cappell-shaneson-codim2}, and was refined in~\cite{COT} to detect more subtle obstructions for knots in $S^3$.  We refer the reader to Carolyn Otto's thesis \cite{otto} for an illustration of why it makes sense to fix a Poincar\'e embedding type {\em before} setting up a surgery problem.  She considers the COT filtration of the string link concordance group, $\mathcal{SL}$, which does not take Poincar\'e embedding type into account, and cleverly shows the non-triviality of $\mathcal{SL}_{n.5}/\mathcal{SL}_{n+1}$ by constructing elements with different Poincar\'e embedding types (i.e. elements distinguished by their $\overline{\mu}$-invariants).  In contrast, there is a unique Poincar\'e embedding type for knots in $S^3$, and it is still an open question if $\mathcal{F}_{n.5}/\mathcal{F}_{n+1}$ is nontrivial for any  $n$ for the COT filtration $\mathcal{F}$ of the knot concordance group.  Our invariants can be recognized as extensions of the homotopy theoretic invariants of Orr for links in $S^3$ \cite{orr2, orr1}, which are equivalent to Milnor's $\overline{\mu}$-invariants.  We note that our invariants are not the first extension of homotop theoretic invariants to knots in non-simply connected manifolds.  In \cite{miller2} Miller extends Milnor's $\overline{\mu}$-invariants to knots in closed, orientable, aspherical, irreducible, Seifert fibered 3-manifolds that are homotopic to the Seifert fiber, although he considers knots up to $\Z \left[ \dfrac{\pi_1(M)}{\Z} \right]$-homology concordance. \\

Fix a knot $J$ in $M$.   The inclusion of the exterior of $J$ into $M$, $\EJ \hookrightarrow M$, induces an epimorphism $\gJ : \pi_J := \pi_1(\EJ) \to G := \pi_1(M)$ with kernel $\G \gJ$.  For $n \geq 2$ let $B_n$ be the mapping cylinder of the induced map $\partial \EJ \to K\left( \dfrac{ \pi_J}{\G_n \gJ}, 1 \right)$, where $\G_n \gJ$ is the $n^{th}$-lower central subgroup of $\G \gJ$ and $K\left( \dfrac{ \pi_J}{\G_n \gJ}, 1 \right)$ is an Eilenberg-MacLane space for $\dfrac{ \pi_J}{\G_n \gJ}$.  Form $X_n$, an approximation of the decomposition of $M$ as $M = \EJ \cup_{\partial \EJ} \overline{N(J)}$, where $\overline{N(J)}$ is a closed tubular neighborhood of $J$, by gluing a solid torus to $B_n$ along $\partial \EJ$ such that $\mu_J$ bounds a disk.  These $\{ X_n \}_{n \geq 2}$ fit into a sequence up to homotopy
$$\xymatrix{ {} & \ar[r] & X_{n+1} \ar[r] & X_n \ar[r] & \cdots \ar[r] & X_2 ~.}$$
Given a knot $K$ in $M$ homotopic to $J$, if there exists a map $\tau: \EK \to B_n$ making the following diagram commute
$$
\xymatrix{
\mu_K \ar@{^(->}[d] \ar[rr] & & \mu_J \ar@{^(->}[d] \\
\partial \EK \ar[rr]^{\cong} \ar[dr]_{\tau} & & \partial \EJ \ar[dl] \\
& B_n
}
$$
where $\mu_J$ and $\mu_K$ are the meridians of $J$ and $K$, respectively, then we may extend $\tau$ to $\ka_n(K, \tau):M \to X_n$.  Our invariant, the homotopy class of $\ka_n(K, \tau)$ in $[M, X_n]$, compares the decompositions $M = \EK \cup_{\partial \EK} \overline{N(K)}$ and $M = \EJ \cup_{\partial \EJ} \overline{N(J)}$ by obstructing the existence of a map $\EK \to B_n$ that extends to $M \to X_n$. \\

We call $\{ \G_n \g_A \}_{n\geq 1}$ the {\em $G$-lower central series of $\g_A$}.  The following Stallings'-type result (see \cite{stallings} for more on Stallings' theorem), well known to experts, implies that concordant knots have isomorphic $G$-lower central series quotients:

\begin{intro isom lcc}
If $f: A \to B$ is a homomorphism over $G$ and $f_{i}: H_{i}(A; \Z[G]) \rightarrow H_{i}(B; \Z[G])$ is an isomorphism for $i=1$ and an epimorphism for $i = 2$ then $f$ induces isomorphisms
$$f: \dfrac{A}{\G_n \g_A} \to \dfrac{B}{\G_n \g_B}$$
for all $n \in \Z$. 
\end{intro isom lcc}

Call a map $\EK \to B_n$ that extends over $M$ as mentioned above {\em extendable}.  We show that if $K$ admits an extendable map to $B_n$ then it induces an isomorphism of the first $n$ $G$-lower central series quotients  in a way that preserves certain peripheral data:

\begin{intro ext. induces isom.}
If $\tau: \EK \to B_n$ is an extendable map for a knot $K_{\beta}$ then $\tau$ induces an isomorphism
$$\tau_*: \dfrac{\PK}{\HKjj} \cong \dfrac{\PJ}{\HJjj} $$
over $G$ with $\tau_* \circ \dagger_K\left( H_K \right)= \dagger_{j+1}(H_J)$ and $\tau_*(\mK) = \mJ$ for all $j < n$.
\end{intro ext. induces isom.}

An {\em $n$-concordance} between $K$ and $J$ is a properly embedded surface in $M\times I$ cobounded by $J$ and $K$ that ``looks'' like an annulus modulo {\em its} $n^{th}$ $\pi_1(M)$-lower central series quotient. 

\begin{intro n-conc iff nth invt}
A knot $K$ is $n$-concordant to $J$ if and only if some $j^{th}$ $\ka$-invariant $\ka_j(K, \tau)$ is defined and trivial for all $j\leq n$.
\end{intro n-conc iff nth invt}

Since a concordance is an $n$-concordance for all $n$, we obtain the following corollary.

\begin{intro ka is conc invt}
If $K$ is concordant to $J$ then for each $n$ there is an extendable map $\tau$ such that $\ka_n(K, \tau)$ is defined and trivial. 
\end{intro ka is conc invt}

Finally, recall that a knot $K$ in $M$ is {\em $J$-characteristic} if there is a continuous degree-one map $\alpha: M \to M$ such that $\alpha(K) = J$ and $\alpha(M-K) \subseteq M-J$.

\begin{intro dense J}
Let $K$ be a knot in a closed, orientable, aspherical 3-manifold.  If $K$ is $J$-characteristic via a map $\alpha:M \to M$ that indices the identity on $G$ then $K$ has trivial $\ka$-invariants.
\end{intro dense J}

As a special case of this, let $L \subset S^3$ be a knot and let $\eta$ be a curve in $\EJ$ that bounds a disk in $M$.  Denote by $J(\eta, L)$ the satellite of $J$ formed by
$$\left( M - N(\eta) \right) \cup -\left( S^3 - N(L) \right), $$
where if $\lambda_{\eta}$ and $\lambda_L$ denote the longitudes of $\eta$ and $L$, respectively, then $\mu_L \sim \lambda_{\eta}^{-1}$ and $\lambda_L \sim \mu_{\eta}$.  This manifold is diffeomorphic to $M$ via a map that is the identity outside a regular neighborhood of the disk bounded by $\eta$.  In particular, this diffeomorphism induces the identity homomorphism on $G$.  

\begin{intro char has triv}
Let $J$ be a homotopically essential knot in a closed, orientable, aspherical 3-manifold.  Let $\eta$ be a curve in $\EJ$ that bounds an embedded disk in $M$ and let $L$ be a knot in $S^3$.  Then $J(\eta, L)$ has trivial $\ka$-invariants.
\end{intro char has triv}

While we expect these invariants to be highly nontrivial for knots in most closed, oriented, prime 3-manifolds, there are a few cases where we expect the invariants to be trivial.  For example, if $M = S^1 \times S^2$ and $K$ generates $\pi_1(S^1 \times S^2) \cong \Z$ or if $M$ is a lens space with prime-order fundamental group and $K$ is homotopically essential then these invariants will likely be trivial.  Indeed, if $H_1\left(\EK; \Z[G] \right) = \Z$ these invariants will be trivial since $\Z$ has trivial lower central series.  However, if $M$ is a lens space with fundamental group of composite order and $K$ does not generate $\pi_1(M)$, if $K$ is a knot in a manifold of the second type, or if $K$ is a knot in $S^1 \times S^2$ representing a class with norm at least three in $\pi_1(S^1 \times S^2) \cong \Z$ then these invariants will likely be highly nontrivial.  The reason is that these invariants essentially measure the homotopy class of the boundary of a tubular neighborhood of the knot in the target space $X_n$, and manifolds with nontrivial fundamental groups allow knots to link themselves nontrivially, creating interesting toral embeddings.  \\

The paper is organized as follows:  we review some basic definitions in section \ref{Knots and concordance:1}, and give the analogous definitions for based spaces in section \ref{Knots and concordance:2}.   We define the category of groups over $G$ and prove Proposition \ref{isom lcc} in section \ref{The $G$-lower central series}.  In section \ref{The invariant} we define the $\ka$-invariants.  We construct the spaces $\{ X_n \}$ in section \ref{Boundary conditions and knot complexes}.  We define extendable maps and prove Proposition \ref{ext. induces isom.} in section \ref{extendable maps sect.}.  We investingate the indeterminacy of the $\ka$-invariants in section \ref{indeterminacy of ka-invariants}.  We introduce $n$-concordance and prove Theorem \ref{n-conc iff nth invt} and Corollary \ref{ka is conc invt} in section \ref{$n$-concordance sect.}.  Finally, in section \ref{Characteristic knots sect.} we consider $J$-characteristic knots in manifolds of the second type and prove Theorems \ref{dense J} and \ref{char has triv}.\\


\noindent {\bf Acknowledgments.}
This work of deducing a coherent model for studying free isotopy classes of circles in general 3-manifolds started several years ago.  I am greatly indebted to my advisor Kent Orr for many helpful conversations and for his overwhelming encouragement.

\section{Preliminaries}
\label{Preliminaries}
Henceforth we work in the smooth category and consider knots in closed, oriented, prime 3-manifolds.  In particular, we allow homologically nontrivial knots.  In Section \ref{Knots and concordance:1} we present some necessary definitions from knot theory.  We give analogues of these same definitions for based knots in Section \ref{Knots and concordance:2}.  In Section \ref{The $G$-lower central series} we define the {\em $G$-lower central series} and prove a Stallings' theorem for the $\M$-lower central series.  Recommended references are \cite{gompf-stip}, \cite{munkres}, and \cite{rolfsen}.

\subsection{Knots and concordance} 
\label{Knots and concordance:1}
$\phantom{M}$

A {\em knot} $K$ in a 3-manifold $M$ is an oriented one-dimensional closed submanifold.  We will assume that $M$ has base point $p \in M$ and that $K$ does not contain $p$.  We also assume that $K$ is based via an embedded path from $p$, although the choice of path will not matter for the results of this paper.  The {\em exterior} of $K$ is the complement of an open normal neighborhood of $K$, $\EK = M-N(K)$, based at $p$.  Note that the inclusion $i_K: \EK \hookrightarrow M$ induces an epimorphism $\gK: \pi_1 (\EK, p) \rightarrow \pi_1 (M, p)$.  We denote the fundamental group of $\EK$ by $\PK$ and the fundamental group of $M$ by $G$. \\

A {\em meridian} $\mK$ of $K$ is an embedded curve in $\partial \EK$ that represents a primitive element of $H_1(\partial \EK)$ and bounds a disk in $N(K)$.  It is uniquely determined up to isotopy in $\partial \EK$.  We will abuse notation by letting $\mK$ denote the curve in $M$, its isotopy class in $\partial \EK$, its homology class in $H_1(\EK)$, and its homotopy class in $\PK$, where in the last case we regard $\mK$ as being based via the path basing $K$. \\

A knot $K$ in $M$ is {\em $L$-characteristic} if there is a continuous degree-one map $\alpha: M \to M$ such that $\alpha(K) = L$ and $\alpha(M-K) \subseteq M-L$. \\

Two knots $K$ and $L$ are said to be {\em cobordant} if there is an oriented submanifold $V$ of $M\times I$ that meets the boundary of $M\times I$ transversely with $V \cap (M\times \{0\}) = L \times\{0\}$ and $V \cap (M\times \{1\}) = K \times\{1\}$, and such that the orientations of $L \times \{0\}$ and $K \times \{1\}$ agree with the orientations induced by $V$.  \\

If $K$ and $L$ are cobordant then we can always choose $V$ to be disjoint from $\{p\} \times I$.  We call $E_V = M\times I - N(V)$ the {\em exterior} of $V$ in $M\times I$, where $N(V)$ is an open normal neighborhood of $V$ that is disjoint from $\{p\} \times I$.  A {\em meridian} $\mu_V$ of $V$ is an embedded curve in $\partial \EV$ that is isotopic to $\mK$.  The path $\rho(t) := p \times t$ induces an isomorphism $\rho_*: \pi_1(\EV, p\times 1) \to \pi_1(\EV, p\times 0)$ by $\rho_*(f) = \rho * f * \rho^{-1}$ (where we read concatenation of paths from left to right).  We therefore write $\pi_1(\EV, p\times 1) = \pi_1(\EV, p\times 0)$ and denote both of these groups by $\PV$. \\

If $V$ is a concordance between $K$ and $L$ then the inclusions $\iota_L : E_L \hookrightarrow E_V$ and $\iota_K : E_K \hookrightarrow E_V$ induce homomorphisms $\iota_L : \pi_L \to \pi_V$ and $\iota_K : \pi_K \to \pi_V$, respectively, and the inclusion $i_V: E_V \hookrightarrow M\times I$ induces an epimorphism $\g_V : \PV \rightarrow G = \pi_1(M)$ such that the following diagram commutes
$$\xymatrix{
\pi_{L} \ar[r]^{\iota_{L}} \ar[rd]_{\g_{L}} & \PV \ar[d]^{\g_V} & \pi_{K} \ar[l]_{\iota_{K}} \ar[ld]^{\g_{K}}\\
& \M.
}
$$

We call a cobordism $V$ a {\em $\M$-homology concordance} if the inclusions $\iota_L : E_L \hookrightarrow E_V$ and $\iota_K : E_K \hookrightarrow E_V$ induce isomorphisms on $H_*( -; \Z[\M])$.  We call $V$ a {\em concordance} if it is homeomorphic to $S^1 \times I$.  A Mayer-Vietoris sequences argument shows that a concordance is a $\M$-homology concordance (see also Proposition 2.1 of  \cite{L2-paper}). \\

\subsection{Based knots and based concordance} 
\label{Knots and concordance:2}

Recall from the previous section that $p \in M$ is the base point of $M$. \\

A {\em based knot} $K_{\alpha}$ is a knot $K$ together with a choice of embedded path $\alpha$ in $M$ from $p$ to $K$ that intersects $K$ only at the endpoint.  We call $\alpha$ a {\em basing} of $K$.  If $K_{\alpha}$ and $K_{\beta}$ denote a knot $K$ based along paths $\alpha$ and $\beta$, respectively, then conjugating by $\alpha \beta^{-1}$ defines a homomorphism $\sigma: \pi_1(M, p) \to \pi_1(M, p)$ such that $\sigma([K_{\alpha}]) = [K_{\beta}]$.  When the basing is clear from the context we will drop the subscript and simply write $K$.  \\

We may assume that a given basing $\alpha$ of $K$ intersects the open normal neighborhood $N(K) \subset M$ in an interval.  Therefore, the {\em exterior} of $K_{\alpha}$ is the exterior $\EK$ of the unbased knot $K$, with base point $p$.  We abuse notation by writing $\partial \EK \cup \alpha$ for the space $\left( \partial \EK \cup \alpha \right) - \left( N(K) \cap \alpha \right)$.  Moreover, whenever we refer to $\partial \EK \cup \alpha$ we assume that this ``truncated" part of alpha is reparameterized as a path with domain $I$ and $\alpha(1) \in \partial \EK$.  The inclusion $\partial \EK \cup \alpha \hookrightarrow \EK$ induces a homomorphism $\pi_1( \partial \EK \cup \alpha, p) \to \PK$. \\

A {\em meridian} of the based knot $K_{\alpha}$ is a loop ${\mK}_{\alpha} = \alpha * \mK * \alpha^{-1}$ in $\partial \EK \cup \alpha$, where $\mK$ is a meridian for the unbased knot $K$ that passes through $\alpha(1)$.  
We will usually drop the subscript $\alpha$ when referring to a meridian or longitude of a based knot. \\

We say that $L_{\alpha}$ and $K_{\beta}$ are {\em based cobordant} (resp. {\em based concordant}) if there is a cobordism (resp. concordance) $V$ between $K$ and $L$ such that $\mu_{K}$ and $\mu_{L}$ have the same image in $\PV$ under $\PK \to \PV$ and $\pi_L \to \PV$, respectively.  If $L_{\alpha}$ and $K_{\beta}$ are based cobordant (resp. based concordant) then they are cobordant (resp. concordant).  However, the converse is not true.  It is often possible to find two different basing $\beta_1$ and $\beta_2$ of $K$ such that $K_{\beta_1}$ is not homotopic relative $p$ to $K_{\beta_2}$.  We say that $L_{\alpha}$ and $K_{\beta}$ are $\M$-homology concordant if they are $\M$-homology concordant as unbased knots.  It can be shown that if $L$ is concordant to $K$ (as unbased knots) via a concordance $C$ then given any basings $\alpha$ and $\beta$ of $L$ and $K$, respectively, there is an inner automorphism $\sigma$ of $\PC$ such that $\sigma(\mK) = \mu_L$. \\

\subsection{The $G$-lower central series}
\label{The $G$-lower central series}
For more on the category $\mathcal{G}^G$ see \cite{L2-paper}.

\begin{defn}
\label{cat of gps over G}
Fix a group $G$.  Let $\mathcal{FG}^G$ denote the category whose objects are group homomorphisms $\g_A: A \to G$ and whose morphisms $f: \g_A \to \g_B$ are group homomorphisms $f: A \to B$ making the following diagram commute:
$$
\xymatrix{
A \ar[rr]^f \ar[dr]_{\g_A} & & B \ar[dl]^{\g_B} \\
& G
}
$$
Let $\mathcal{G}^G$ denote the full subcategory of $\mathcal{FG}^G$ whose objects are epimorphisms $\g_A: A \to G$.  We call morphisms in $\mathcal{FG}^G$ and $\mathcal{G}^G$ {\em homomorphisms over $G$}.
\end{defn}

\begin{defn}
For an object $\g_A \in \mathcal{FG}^G$, define
$$\G \g_A = \text{Ker}\{ \g_A \}.$$
\end{defn}

\begin{defn}
\label{def of lower cent series}
Define the {\em $G$-lower central series of an object $\g_A \in \mathcal{FG}^G$} by 
$$\G_1 \g_A = \G \g_A $$
and 
$$\G_{n+1} \g_A = [\G_1 \g_A, \G_n \g_A].$$
\end{defn}

In particular, the $G$-lower central series of $\g_A$ is precisely the lower central series of $\G \g_A$.  For convenience, we will drop the subscript of $\G_1$ and simply write $\G \g_A$ for $\G_1 \g_A$.  It is clear that for $\g_A: A \to G$ this defines a normal series
$$\xymatrix{ A \ar@{}[r]|{\unrhd} & \G \g_A  \ar@{}[r]|{\unrhd} & \G_2 \g_A  \ar@{}[r]|{\unrhd} & \cdots  \ar@{}[r]|{\unrhd} & \G_n \g_A  \ar@{}[r]|{\unrhd} & \cdots . }$$

\begin{prop}
\label{isom lcc}
If $f: \g_A \to \g_B$ is a morphism in $\mathcal{G}^G$ and $f_{i}: H_{i}(A; \Z[G]) \rightarrow H_{i}(B; \Z[G])$ is an isomorphism for $i=1$ and an epimorphism for $i = 2$ then $f$ induces isomorphisms
$$f: \dfrac{A}{\G_n \g_A} \to \dfrac{B}{\G_n \g_B}$$
for all $n \in \Z$. 
\end{prop}
\begin{proof}
Because $H_i(\G \g_A) = H_{i}(A; \Z[G])$ and $H_i(\G \g_B) = H_{i}(B; \Z[G])$, Stallings' Theorem implies that $f$ induces an isomorphism $\dfrac{\G \g_A}{\G_n \g_A} \to \dfrac{\G \g_B}{\G_n \g_B}$ for all $n$.  The result now follows by applying the Five Lemma to the following commutative diagram of short exact sequences induced by $f$:
$$
\xymatrix{
0 \ar[r] & \dfrac{\G \g_A}{\G_n \g_A} \ar[r] \ar[d]^{\cong} & \dfrac{A}{\G_n \g_A} \ar[r]^{\g_A} \ar[d] & G \ar[r] \ar[d]^= & 0 \\
0 \ar[r] & \dfrac{\G \g_B}{\G_n \g_B} \ar[r] & \dfrac{B}{\G_n \g_B} \ar[r]^{\g_B} & G \ar[r] & 0
}
$$
\end{proof}

\begin{example}
\label{H ex}
Consider the following short exact sequence
$$\xymatrix{0\ar[r] & F_2 \ar[r]^(.4){\iota} & F\times \Z \ar[r]^(.5){\g} & \dfrac{F}{F_3} \ar[r] & 0 }$$
where $F$ is the free group on two generators $x$ and $y$, $\Z$ is generated by $t$, and $F_2$ is the commutator subgroup of $F$.  Define $\iota: F_2 \rightarrow F\times \Z$ by $\iota([x, y]^{\omega}) = [x, y]^{\omega}t^{-1}$ for any word $\omega \in F\times \Z$.  Then $\G \g \cong F_2$ is the free group on infinitely many generators and $\G_n \g \cong (F_2)_n$.  One can show that $\dfrac{F}{F_3}$ is the fundamental group of the Heisenberg manifold, a circle bundle over the torus, and that $F \times \Z$ is the fundamental group of the complement of a fiber in this manifold. 
\end{example}


\section{Homotopy invariants}
Fix a based knot $J_{\alpha}$ in $M$; this is the knot to which we will compare all other based knots.  In Section \ref{Boundary conditions and knot complexes} we will explicitly construct an infinite tower of spaces
$$\xymatrix{ \ar[r] & X_n \ar[r] & X_{n-1} \ar[r] & \cdots \ar[r] & X_2 }$$
{\em that depends on $J_{\alpha}$}.  The $X_n$ can be thought of as successive approximations of $\EJ$.  Our goal in Section \ref{The invariant} is to associate to any suitable unbased knot $K$ a homotopy class $\ka \in [M, X_n]$.  To do this, we first associate to any suitable based knot $K_{\beta}$ a homotopy class $\ka \in [M, X_n]_0$.  We discuss in \ref{indeterminacy of ka-invariants} the extent to which $\ka$ depends on the basing $\beta$. The definition of these invariants requires the notion of an {\em extendable map}, which, for brevity, we put off defining until Section \ref{extendable maps sect.}.  Roughly speaking, two knots are {\em $n$-concordant} if they cobound a surface $V$ that ``looks'' like an annulus modulo $\G_n \g_V$. We investigate the intimate relationship between $\ka$-invariants and $n$-concordance in Section \ref{$n$-concordance sect.}.  Finally, we investigate the relationship between $\ka$-invariants and $J$-characteristicness in Section \ref{Characteristic knots sect.}.

\subsection{The invariant}
\label{The invariant}

\begin{defn}
A {\em boundary condition} is a pair $(\dagger, \mu)$ consisting of a homomorphism $\dagger: H \to \pi$, where $H$ is a free abelian group on two elements, and a primitive element $\mu \in H$.  Note that the boundary condition does not carry with it a particular choice of isomorphism $H \cong \Z^2$, and $\dagger$ need not be injective.  We abuse notation by simply writing $\mu$ for $\dagger(\mu)$.  Later, we will have a commutative diagram,
$$
\xymatrix{
H \ar[rr]^{\dagger} \ar[dr] && \pi \ar[dl] \\
& G
}
$$
making $\dagger$ a morphism in $\mathcal{FG}^G$.
\end{defn}

\begin{defn}[Whitehead \cite{whitehead}, p. 244]
Given a group $\pi$, an {\em Eilenberg-MacLane space of type $(\pi, 1)$} is a based space $K(\pi, 1)$ such that $\pi_n \left( K(\pi,1) \right) = 0$ for all $n>1$ together with an isomorphism $\pi_1 \left( K(\pi,1) \right) \to \pi$.
\end{defn}

\begin{remark}
\label{construct X}
We now associate to the boundary condition $(\dagger, \mu)$ a space that is well-defined up to based homotopy equivalence.
Let $T^2$, a torus, be an Eilenberg-MacLane space for $H$, and let $*$ denote the base point of $T^2$.  With some abuse of notation, let $\mu$ be an embedded curve in $T^2$ containing $*$ that represents the class $\mu \in H$.  Let $K(\pi, 1)$ be an Eilenberg-MacLane space for $\pi$, based at a point $b$.  The homomorphism $\dagger$ determines a based map $\dagger: T^2 \to K(\pi, 1)$, up to based homotopy.  
Let $B_{\dagger}$ be the mapping cylinder of $\dagger$,
$$ B_{\dagger} = \dfrac{T^2 \times I \coprod K(\pi, 1)}{(x, 0) \sim \dagger(x)} \phantom{a}. $$
We will abuse notation by denoting $T^2 \times \{1\}$ and $\mu \times \{1\}$ by $T$ and $\mu$, respectively.  Define
$$X_{\dagger} = B_{\dagger} \cup ST, $$
where the solid torus $ST$ is attached to $T$ via a homeomorphism $\partial (ST) \cong T$ so that $\mu$ bounds a disk in $ST$.  Both $B_{\dagger}$ and $X_{\dagger}$ are based at the point $b \in K(\pi, 1)$.  We connect $T$ and $\mu$ to the base point $b$ via the straight-line path along $T^2 \times I$ from $(\ast, 0) \sim b$ to $(\ast, 1)$.  Since homotopy pushouts are homotopy invariant, the space $X_{\dagger}$ is well defined up to based homotopy type.  Hence, we drop the subscript and simple write $X$ for $X_{\dagger}$. 
\end{remark}

\begin{defn}
We call $X$ the {\em complex of the boundary condition $(\dagger, \mu)$}.  
\end{defn}

\begin{defn}
\label{x_n}
Let $K_{\beta}$ be a based knot.  Let $H_K = \pi_1(\partial \EK \cup \beta, p)$ and let $\dagger_K: H_K \to \pi_K$ be the homomorphism induced by $\partial \EK \cup \beta \hookrightarrow \EK$.  Call $\left( \dagger_K, \mu_K \right)$ the {\em boundary condition associated to $K_{\beta}$}.  $\dagger_K$ may be viewed as a morphism in $\mathcal{FG}^G$.
\end{defn}

\begin{defn}
Fix a based knot $J_{\alpha}$; this knot will remain fixed for the duration of the article.  For $n \geq 2$, define
$$ \dagger_n: H_J \to \PJ \to \dfrac{\PJ}{\HJn}$$ 
to be $\dagger_J$ followed by the canonical quotient homomorphism.  Define the {\em $n^{th}$-knot complex of $J_{\alpha}$}, denoted $X_n$, to be the knot complex associated to the boundary condition $\left( \dagger_n, \mJ \right)$.  Denote the mapping cylinder of $T^2 \to K\left( \dfrac{\PJ}{\HJn}, 1 \right)$, as in Remark \ref{construct X}, by $B_n$.  Let $T_n$ and $\mu_J$ denote $T^2 \times \{1\}$ and $\mJ \times \{1\}$, respectively.  Denote the base point of $B_n$ by $b$, for all $n$.  We abuse notation by letting $\gJ$ denote the epimorphism of fundamental groups $\dfrac{\PJ}{\HJn} \to \pi_1(X_n)$ induced by the inclusion $B_n \to X_n$.
\end{defn}

The spaces $\{X_n\}_{n \geq 2}$ were defined with respect to the fixed knot $J$ and can be thought of as successive approximations of the decomposition $M = \EJ \cup_{\partial N(J)} \overline{N(J)}$.  Given a knot $K$ (freely) homotopic to $J$ in $M$, the invariant below compares $K$ to $J$ by comparing the decomposition $M = \EK \cup_{\partial N(K)} \overline{N(K)}$ to $X_n$.  

\begin{defn}
\label{ka invts}
Let $K_{\beta}$ be a based knot such that $[K_{\beta}]$ is conjugate to $[J_{\alpha}]$ in $G$.  Let $\tau: E_K \rightarrow B_n$ be an extendable map for $K_{\beta}$ as in Definition \ref{extendable map}, with epimorphism $\gK: \PK \to G$.  Then it extends to a based map $\ka_n(K,\tau) : M \rightarrow X_n$ by Proposition \ref{ext extends}, and we obtain a based homotopy class $[\ka_n(K,\tau)] \in [M, X_n]_0$.  We refer to the based class $[\ka_n(K,\tau)] \in [M, X_n]_0$ as an {\em $n^{th}$ $\ka$-invariant for the based knot $K_{\beta}$}.  We define the free homotopy class
$$[\ka_n(K,\tau)] \in [M, X_n]$$
to be the {\em $n^{th}$ $\ka$-invariant of the pair $(K, \tau)$}.  Let $\tau_n: \EJ \to B_n$ be as in Example \ref{can. extension}.  We say that $[\ka_n(K,\tau)]$ is {\em $J$-trivial} (or simply {\em trivial} when $J$ is clear from the context) if $[\ka_n(K,\tau)] = [\ka_n(J,\tau_n)]$ in $[M, X_n]$.  That is, if $\ka_n(K,\tau)$ is freely homotopic to $\ka_n(J,\tau_n)$. 
\end{defn}

\subsection{Boundary conditions and knot complexes}
\label{Boundary conditions and knot complexes}
We now construct the tower of spaces
$$\xymatrix{ \ar[r] & X_n \ar[r] & X_{n-1} \ar[r] & \cdots \ar[r] & X_2 ~.}$$

\begin{defn}
Let $(\dagger_0:H_0 \to \pi_0, \mu_0)$ and $(\dagger_1: H_1 \to \pi_1, \mu_1)$ be boundary conditions.  A {\em homomorphism (resp. isomorphism) of boundary conditions} 
$$p: (\dagger_0, \mu_0) \to (\dagger_1, \mu_1)$$
is a homomorphism (resp. isomorphism) $p: \pi_0 \to \pi_1$ with $p(\mu_0) = \mu_1$ and which restricts to an isomorphism $p: \dagger_0(H_0) \cong \dagger_1(H_1)$.
\end{defn}

\begin{remark}
Since $H_0$ and $H_1$ are both isomorphic to $\Z^2$, given any homomorphism $p: \pi_0 \to \pi_1$ that restricts to an isomorphism $\dagger_0(H_0) \cong \dagger_1(H_1)$ with $p(\mu_0) = \mu_1$, there exists an isomorphism $H_0 \to H_1$ making the diagram
$$
\xymatrix{
H_0 \ar@{..>}[r]^{\cong} \ar[d]_{\dagger_0} & H_1 \ar[d]^{\dagger_1}\\
\pi_0 \ar[r]^p &  \pi_1
}
$$
commute.  We will freely make use of this fact in the future without explicitly choosing such an isomorphism.
\end{remark}

\begin{prop}
\label{map of class space is WD}
Let $p: (\dagger_0, \mu_0) \to (\dagger_1, \mu_1)$ be a homomorphism of boundary conditions, let $X_i$ be the complex of $(\dagger_i, \mu_i)$, and let $T_i$ and $\mu_i$ denote $T$ and $\mu$ in $X_i$, respectively (see Remark \ref{construct X}).  Then $p$ induces a based map $P: X_0 \to X_1$, well defined up to based homotopy, such that $P\left( K(\pi_0, 1) \right) \subset K(\pi_1, 1)$ and $P(T_0) = T_1$ with $P(\mu_0) = \mu_1$.  Moreover, if $p$ is an isomorphism then $P$ is a homotopy equivalence.
\end{prop}
\begin{proof} 
For $i \in \{0, 1 \}$, let $\dagger_i: T^2 \to K(\pi_i, 1)$ be a based map induced by $\dagger_i$, let $\mu_i \subset T^2$ be an embedded curve representing $\mu_i \in H_i$, and let $X_{\dagger_i} = B_{\dagger_i} \cup ST$, as above.  The homomorphism $p: \pi_0 \to \pi_1$ induces a map $Q: K(\pi_0, 1) \to K(\pi_1, 1)$ up to based homotopy, and hence a map $Q: X_{\dagger_0} \to X_{Q \circ \dagger_0}$, up to based homotopy, with $Q(T_0) = T_{Q \circ \dagger_0}$ and $Q(\mu_0) = \mu_{Q \circ \dagger_0}$, and which is easily seen to be a homotopy equivalence if $p$ is an isomorphism.  Since $Q_* = p: \dagger_0(H_0) \to \dagger_1(H_1)$ isomorphically with $p(\mu_0) = \mu_1$, there is a based homotopy equivalence $R: X_{Q \circ \dagger_0} \to X_1$ that maps $K(\pi_1, 1)$ to itself and $T_{Q \circ \dagger_0}$ to $T_1$ with $\mu_{Q \circ \dagger_0} \mapsto \mu_1$.  We define $P:= R \circ Q.$ 
\end{proof}

\begin{lem}
\label{X_n = G}
For each $n$, there is a canonical isomorphism $\pi_1(X_n, b) \cong G$.
\end{lem}
\begin{proof}
Recall that $X_n = B_n \cup_{T_n} ST$, where $ST$ is attached so that $\mJ$ bounds a disk.  The Seifert-van Kampen Theorem implies that $\pi_1(X_n, b)$ is the pushout of the following diagram,
$$
\xymatrix{
H_J \ar[r]^{\dagger_n} \ar[d] & \dfrac{\PJ}{\HJn} \\
\Z
}
$$
where the kernel of the vertical arrow is the subgroup generated by $\mJ$.  Hence, there is a unique homomorphism $\pi_1(X_n, b) \to G$ making the following diagram commute,
$$
\xymatrix{
H_J \ar[r]^{\dagger_n} \ar[d] & \dfrac{\PJ}{\HJn} \ar[d] \ar@/^/[ddr] \\
\Z \ar[r] \ar@/_/[drr] & \pi_1(X_n, b) \ar@{..>}[dr] \\
& & G
}
$$
where the bottom-left homomorphism $\Z \to G$ sends the generator of $\Z$ to the homotopy class of $J_{\alpha}$ and $\dfrac{\PJ}{\HJn} \to G$ is the canonical epimorphism induced by $\gJ: \pi_J \to G$.  Since $\mJ$ normally generates $\dfrac{\HJ}{\HJn}$ in $\dfrac{\PJ}{\HJn}$, $\pi_1(X_n, b) \to G$ is an isomorphism. 
\end{proof}

\begin{remark}
Henceforth we fix this isomorphism between $\pi_1(X_n, b)$ and $G$, and write $\pi_1(X_n, b) = G$.  
\end{remark}

\begin{remark}
\label{canonical J}
The canonical epimorphisms $q_{n+1, n}: \dfrac{\PJ}{\HJnn} \to \dfrac{\PJ}{\HJn}$ induce a tower of boundary conditions
$$\xymatrix{ \ar[r] & \left( \dagger_{n+1}, \mJ \right) \ar[r] & \left( \dagger_n, \mJ \right) \ar[r] & \cdots \ar[r] & \left( \ \dagger_2, \mJ \right). }$$ 
By Proposition \ref{map of class space is WD}, this tower of boundary conditions induces a tower of knot complexes
$$\xymatrix{ \ar[r] & X_{n+1} \ar[r]^{Q_{n+1,n}} & X_n \ar[r] & \cdots \ar[r] & X_2~. }$$ 
Each map in the tower induces the identity homomorphism on $G$ and is defined up to based homotopy equivalence.  We therefore obtain maps of sets
$$ Q_{n+1,n}: [M, X_{n+1}]_0 \to [M, X_n]_0 $$ 
defined by $Q_{n+1,n}\left( [f] \right) := [Q_{n+1,n} \circ f]$.
\end{remark}

We place the following theorem here for completeness.  Extendable maps are defined in Definition \ref{extendable map}.
\begin{thrm}
\label{ka down tower}
Let $K_{\beta}$ be a based knot and suppose there is an extendable map $\tau: \EK \to B_{n+1}$ for $K_{\beta}$.  Then there is an extendable map $\tilde{\tau}: \EK \to B_n$ for $K_{\beta}$ such that $Q_{n+1,n}\left( [\ka_{n+1}(K, \tau) ] \right) = [\ka_n(K, \tilde{\tau}) ]$.
\end{thrm}
\begin{proof}
Since $T_{n+1}$ is a sub-complex of the CW-complex $B_{n+1}$ we may assume that $q_{n+1, n}: \left( \dag_{n+1}, \mJ \right) \to \left( \dag_n, \mJ \right)$ induces a map $\tilde{Q}: B_{n+1} \to B_n$ that restricts to a homeomorphism from $T_{n+1}$ to $T_n$ with $\tilde{Q}(\mJ) = \mJ$ and such that $\tilde{Q}\left( B_{n+1} - T_{n+1} \right) \subseteq B_n - T_n$.  Then $\tilde{\tau} := \tilde{Q} \circ \tau: \EK \to B_n$ is the desired extendable map. 
\end{proof}
$\phantom{}$ \\

\subsection{extendable maps}
\label{extendable maps sect.}
We now define extendable maps and show that if $\EK \to B_n$ is extendable then it induces isomorphisms on the first $n$ $G$-lower central series quotients, preserving peripheral data.

\begin{defn}
\label{extendable map}
Let $N$ be a manifold, possibly with boundary, based at $q \in N$, together with an epimorphism $\g_N: \pi_1(N, q) \twoheadrightarrow G$.  Let $V \subset N$ be a codimension-two connected submanifold, not containing $q$, with trivial normal bundle and $\partial V \subset \partial N$, meeting transversely.  Let $V_{\beta}$ denote $V$ with a choice of basing $\beta$.  Let $N(V)$ be an open regular neighborhood of $V$ that intersects $\beta$ in a subinterval, let $\mu_V$ be the meridian of $V$, let $\EV = N - N(V)$ based at $q$, and let $\PV = \pi_1(E_V, q)$.  Let $\g_V: \PV \twoheadrightarrow G$ be the epimorphism of fundamental groups induced by the inclusion of $\EV$ into $N$ followed by $\g_N$.  We call a based map $f: \EV \to B_n$ {\em extendable} if all of the following hold:
\begin{itemize}
\item[(i)] $f(\partial \EV) \subseteq T_n$ and $f(\mu_V) = \mu_J$, 
\item[(ii)] $f_*\left( \pi_1(\partial \EV \cup \beta, q)\right) = \dagger_n(H_J)$, and 
\item[(iii)] $f_*: \PV \to \dfrac{\pi_J}{\G_n \g_{J}}$ is a morphism of $\mathcal{G}^G$, as in the following diagram:
\end{itemize}
$$
\xymatrix{
\PV \ar[rr]^(.55){f_*} \ar[dr]_{\g_V} & & \dfrac{\pi_J}{\G_n \g_J} \ar[dl]^{\g_{J}} \\
& G 
}
$$ 
If $V$ is not connected, we call a based map $f: \EV \to B_n$ extendable if it satisfies (i) - (iii) for each connected component of $V$.  
\end{defn}

\begin{remark} \hspace{.5cm}
\begin{enumerate}
\item The definition of extendable requires both a choice of basing for $V$ and a choice of epimorphism $\g_N: \pi_1(N, q) \to G$.
\item If $V$ is not connected then a choice of basing is required for each connected component of $V$.
\item If $f$ only satisfies (ii) and (iii) of Definition \ref{extendable map} with $f_*(\mu_V) = \mu_J$ then $f$ is based homotopic to an extendable map.
\end{enumerate}
\end{remark}

\begin{prop}
\label{ext extends}
If $f: \EV \to B$ is extendable then it extends to a based map $f: N \to X$.
\end{prop}
\begin{proof} 
We assume that $V$ is connected because it is sufficient to prove the result for each connected component of $V$.  Choose a homeomorphism $\overline{N(V)} \cong V \times D^2$, where $\overline{N(V)}$ is the closure of $N(V)$, such that $\{*\} \times S^1$ is identified with $\mu_V$ for some point $* \in V$.  To prove the result we show that $f: V\times S^1 \to T_n$ extends to a map $F: V \times D^2 \to ST$.  Let $e^0 = *$ be the 0-cell of $V$ and let $\{ e^k_j \}_{j=1}^{n_k}$ be the $k$-cells of $V$ for $1 \le k \le n$.  Let $f^0$ and $f^1$ be the 0-cell and 1-cell of $S^1$, respectively, and let $f^2$ be the 2-cell of $D^2$.  Then $V \times S^1$ has a CW-decomposition with 0-cell $e^0 \times f^0$, 1-cells $\{ e_j^1 \times f^0 \}_{j=1}^{n_1}$ and $e^0 \times f^1$, and 2-cells $\{e_j^2 \times f^0 \}_{j=1}^{n_2}$ and $\{ e_j^1 \times f^1 \}_{j=1}^{n_1}$.  $V \times D^2$ has the same 0- and 1-cells, but one-additional 2-cell, $e^0 \times f^2$, that is attached by a degree-one map to $e^0 \times f^1$.  Since $f(e^0 \times f^1) = \mu$, $F$ extends to a map from $V \times S^1 \cup_{e^0 \times f^1} D^2$ to $ST$.  $F$ extends $f$ over the 2-skeleton of $V \times D^2$.  We may now extend $F$ over the remainder of $V \times D^2$ because $\pi_{k}(ST) = 0$ for all $k \geq 2$.
\end{proof}

\begin{example}
\label{can. extension}
The canonical epimorphism $q_n: \PJ \to \dfrac{\PJ}{\G_n \gJ}$ induces an extendable map $\tau_n: \EJ \to B_n$, well defined up to based homotopy.  In this case the diagram in (iii) of Definition \ref{extendable map} is
$$
\xymatrix{
\PJ \ar[rr]^{q_n} \ar[dr]_{\g_J} & & \dfrac{\PJ}{\HJn} \ar[dl]^{\gJ} \\
& G 
}
$$ 
\end{example}

\begin{prop}
\label{basics of \EK}
Let $K$ be a knot in a closed, oriented, prime 3-manifold $M$.
\begin{enumerate}
%
\item If either $M$ is irreducible or $M \cong S^1 \times S^2$ and $[K] = 0 \in \pi_1(S^1 \times S^2) \cong \Z$ then the inclusion $\partial \EK \hookrightarrow M = \EK \cup_{\partial \EK} \overline{N(K)}$ induces an isomorphism
$$\xymatrix{ H_1 \left( \partial \EK; \Z[\M] \right) \cong H_1 \left( \EK; \Z[\M] \right) \oplus H_1 \left( \overline{N(K)}; \Z[\M] \right).}$$

\item If $[K]$ has infinite order in $\pi_1(M)$ then 
$$H_2 \left( \EK; \Z[G] \right) = H_2 \left( \Gamma \g_K \right) = 0.$$

\item If $M$ is irreducible and $[K]$ has infinite order in $\pi_1(M)$ then
$$\xymatrix{ H_1 \left( \EK; \Z[\M] \right) \cong H_1 \left( \partial \EK; \Z[\M] \right) \cong H_1\left( \mK; \Z\left[ \frac{\M}{[K]\Z} \right] \right) = \bigoplus_{\frac{\M}{[K]\Z}} \Z }$$ 
as $\Z[\M]$-modules, where $[K]\Z \subset \M$ is the infinite cyclic subgroup generated by $[K]$ and $\dfrac{\M}{[K]\Z}$ is the set of right cosets of $[K]\Z$ in $\M$. 

\item For any knot $K$ in $M$, if $\tau: E_K \to B_n$ is extendible then the image of the homomorphism 
$$H_2 \left( \G \gK \right) \to H_2 \left( \dfrac{\G \gJ}{\G_j \gJ} \right)$$
induced by $\tau$ equals the image of
$$H_2 \left( \G \gJ \right) \to H_2 \left( \dfrac{\G \gJ}{\G_j \gJ} \right)$$
induced by $q_j: \PJ \to \dfrac{\PJ}{\G_j \gJ}$ for all $j \leq n$.
\end{enumerate}
\end{prop}
\begin{proof}
For the first statement consider the following long exact sequence induced from the decomposition $M = \EK \cup_{T^2} \overline{N(K)}$:
$$
\xymatrix{
H_2 \left( \partial \EK; \Z[\M] \right) \ar[r] & \txt{ $H_2 \left( \EK; \Z[\M] \right)$\\ $\oplus$\\ $H_2 \left( \overline{N(K)}; \Z[\M] \right)$} \ar[r] &H_2 \left( M; \Z[\M] \right) &{}\save[]+<2cm,-1cm> *{\text{ $(\bigstar)$}} \restore \\
{} \ar[r] &  H_1 \left( \partial \EK; \Z[\M] \right) \ar[r] & \txt{ $H_1 \left( \EK; \Z[\M] \right)$\\ $\oplus$\\ $H_1 \left( \overline{N(K)}; \Z[\M] \right)$} \ar[r] &H_1 \left( M; \Z[\M] \right) 
}
$$
If $M$ is irreducible then $H_1 \left( M; \Z[G] \right) = H_2 \left( M; \Z[G] \right) = 0$ (being the homology of the universal cover of $M$), establishing (1).  If $M \cong S^1 \times S^2$ then $H_1 \left( S^1 \times S^2; \Z[\Z] \right) = 0$, so we need only consider
$$H_2 \left( S^1 \times S^2; \Z[\Z] \right) \to H_1 \left( \partial \EK; \Z[\Z] \right).$$ 
Since $[K]=0$ the image of $\pi_1( \partial \EK) \to \pi_1(S^1 \times S^2) \cong \Z$ is trivial, so $\partial \EK$ lifts to an infinite number of tori (one for each element of $\Z$) and, if $S$ is a sphere representing a generator of $H_2 \left( S^1 \times S^2; \Z[\Z] \right) \cong H_2 \left( \R \times S^2 \right) \cong \Z$ then $H_2 \left( S^1 \times S^2; \Z[\Z] \right) \to H_1 \left( \partial \EK; \Z[\Z] \right)$ is given by the intersection of $S$ with the lifts of $\partial \EK$.  Since $K$ is null-homotopic in $S^1 \times S^2$ these intersection curves must pair up so as to sum to zero in $H_1 \left( \partial \EK; \Z[\Z] \right)$, making $H_2 \left( S^1 \times S^2; \Z[\Z] \right) \to H_1 \left( \partial \EK; \Z[\Z] \right)$ the zero map.  This establishes (1). \\

If $[K]$ has infinite order in $G$ then $H_2 \left(\partial \EK; \Z[G] \right)=0$.  Moreover, $H_2 \left( \EK; \Z[G] \right) \to H_2 \left( \Gamma \g_K \right)$ is surjective by Hopf's theorem.  If $M$ is irreducible then (2) follows.  If $M \cong S^1 \times S^2$ then $H_1 \left(\partial \EK; \Z[\Z] \right) \cong \Z\left[ \Z_{[K]} \right]$.  By a similar argument to that of (1), if $S$ is a sphere representing a generator of $H_2(\R \times S^2)$ and $t$ generates $\Z_{[K]}$ then
$$H_2( \R \times S^2) \hookrightarrow H_1 \left(\partial \EK; \Z[\Z] \right)$$
by $[S] \mapsto 1+t+\cdots +t^{[K]-1}$,
establishing (2).\\

If $M$ is irreducible and $[K]$ has infinite order in $G$ then $H_1 \left( \overline{N(K)}; \Z[\M] \right) = 0$ and
$$H_1 \left( \partial \EK; \Z[\M] \right) \cong H_1\left( \mK; \Z\left[ \frac{\M}{[K]\Z} \right] \right) = \bigoplus_{\frac{\M}{[K]\Z}} \Z, $$
establishing (3).\\

If $[K]$ has infinite order in $G$ then (4) follows from (2).  Consider the commutative diagram,
$$
\xymatrix{
H_2 \left( \partial \EK; \Z[G] \right) \ar@{->>}[r] \ar[d]^{\cong} & H_2 \left( \EK; \Z[G] \right) \ar@{->>}[r] \ar[d] & H_2 \left( \G \gK \right) \ar[r] \ar[d] & H_2 \left( \dfrac{\G \gK}{\G_j \gK}  \right) \ar[d]\\
H_2 \left( T_n; \Z[G] \right) \ar[r] & H_2 \left( B_n; \Z[G] \right) \ar[r] & H_2 \left( \dfrac{\G \gJ}{\G_n \gJ}  \right) \ar[r] & H_2 \left( \dfrac{\G \gJ}{\G_j \gJ}  \right) \\
H_2 \left( \partial \EK; \Z[G] \right) \ar@{->>}[r] \ar[u]_{\cong} & H_2 \left( \EK; \Z[G] \right) \ar@{->>}[r] \ar[u] & H_2 \left( \G \gK \right) \ar[r] \ar[u] & H_2 \left( \dfrac{\G \gK}{\G_j \gK}  \right) \ar[u]
}
$$
where $H_2 \left( \EK; \Z[G] \right) \to H_2 \left( \G \gK \right)$ is surjective by Hopf's theorem.  If $M$ is irreducible then $H_2 \left( \partial \EK; \Z[G] \right) \to H_2 \left( \EK; \Z[G] \right)$ is surjective by the long exact sequence $(\star)$, establishing (4).  Suppose that $M \cong S^1 \times S^2$ and $K$ is null-homotopic in $M$.  Then from the long exact sequence $(\star)$ we obtain a short exact sequence as follows:
$$
\xymatrix{
0 \ar[r] & H_2 \left( \partial \EK; \Z[\Z] \right) \ar[r] & H_2 \left( \EK; \Z[\Z] \right) \ar[r] & H_2 \left( S^1 \times S^2; \Z[\Z] \right) \ar[r] & 0\\
}
$$
Since $H_2 \left( S^1 \times S^2; \Z[\Z] \right) \cong \Z$ is a free abelian group $H_2 \left( \EK; \Z[\Z] \right)$ splits as a direct sum, 
$$H_2 \left( \EK; \Z[\Z] \right) \cong H_2 \left( \partial \EK; \Z[\Z] \right) \oplus H_2 \left( S^1 \times S^2; \Z[\Z] \right).$$
Hence, we obtain the following commutative diagram:
$$
\xymatrix{
H_2 \left( \partial \EK; \Z[\Z] \right) \oplus H_2 \left( \R^1 \times S^2 \right) \ar[r]^(.65){\cong} \ar@{->>}[d] & H_2 \left( \EK; \Z[\Z] \right) \ar@{->>}[d]\\
H_2 \left( \partial \EK; \Z[\Z] \right) \oplus H_2 \left( \pi_1 (\R^1 \times S^2) \right) \ar[r] & H_2 \left( \PK; \Z[\Z] \right)
}
$$
As $H_2 \left( \pi_1 (\R^1 \times S^2) \right) = 0$, $H_2 \left( \partial \EK; \Z[\Z] \right) \to H_2 \left( \G \g_K \right)$ is surjective and (4) follows.
\end{proof}

\begin{prop}
\label{ext. induces isom.}
If $\tau: \EK \to B_n$ is an extendable map for a knot $K_{\beta}$ then $\tau$ induces an isomorphism
$$\tau_*: \dfrac{\PK}{\HKjj} \cong \dfrac{\PJ}{\HJjj} $$
over $G$ with $\tau_* \circ \dagger_K\left( H_K \right)= \dagger_{j+1}(H_J)$ and $\tau_*(\mK) = \mJ$ for all $j < n$.
\end{prop}
\begin{proof}
It follows from the definition of extendable map that $\tau_*\circ \dagger_K (H_K) = \dagger_{j+1}(H_J)$ for $j<n$ and $\tau_*(\mK) = \mJ$. \\

It remains to show that $\tau_*$ induces an isomorphism between $\dfrac{\PK}{\HKjj}$ and $\dfrac{\PJ}{\HJjj}$.  It is clear that $\tau$ extends to a map $\tau: M \to X_n$ such that $\tau: \partial \EK \to T_n$ homemorphically with $\tau(\mu_K) = \mu_J$, and, since $\tau$ is extendable,
$$
\xymatrix{
H_K \ar[r]^{id} \ar[d]_{\dagger_K} & H_K \ar[d]^{\dagger_K} \ar[r] & H_J \ar[d]^{\dagger_{j+1}} \\
\PK \ar[r] \ar[dr]_{\gK} & \dfrac{\PK}{\HKjj} \ar[r] \ar[d]^{\gK} & \dfrac{\PJ}{\HJjj} \ar[dl]^{\gJ} \\
& G
}
$$
commutes for $j<n$.  Hence, we obtain the following commutative diagram of long exact sequences:
$$
\xymatrix{
{} \ar[r] & H_2 \left( M; \Z[G] \right) \ar[r] \ar[d] & H_2 \left( M, \EK; \Z[G] \right) \ar[r] \ar[d] & H_1 \left( \EK; \Z[G] \right) \ar[r] \ar[d] & 0 \\
{} \ar[r] & H_2 \left( X_n; \Z[G] \right) \ar[r] & H_2 \left( X_n, B_n; \Z[G] \right) \ar[r] & H_1 \left( B_n; \Z[G] \right) \ar[r] & 0\\
{} \ar[r] & H_2 \left( M; \Z[G] \right) \ar[r] \ar[u] & H_2 \left( M, \EJ; \Z[G] \right) \ar[r] \ar[u] & H_1 \left( \EJ; \Z[G] \right) \ar[r] \ar[u] & 0 
}
$$
where $H_1\left(M; \Z[G] \right) = H_1\left(X_n; \Z[G] \right) = 0$, being the homology of the universal covers of $M$ and $X_n$, respectively.  By excision, for $L = K,J$
$$H_2 \left( M, E_L; \Z[G] \right) \cong H_2 \left( N(L), \partial N(L); \Z[G] \right) \cong H_2 \left( ST, T_n; \Z[G] \right) \cong H_2 \left( X_n, B_n; \Z[G] \right).$$
Moreover, the image of
$$H_2 \left( M; \Z[G] \right) \to H_2 \left( M, E_K; \Z[G] \right) \cong H_2 \left( X_n, B_n; \Z[G] \right)$$
equals that of
$$H_2 \left( M; \Z[G] \right) \to H_2 \left( M, E_J; \Z[G] \right) \cong H_2 \left( X_n, B_n; \Z[G] \right).$$
Since
$$H_1 \left( \EJ; \Z[G] \right) \to H_1 \left( B_n; \Z[G] \right)$$
is an isomorphism (recall, $H_1 \left( \EJ; \Z[G] \right)= \dfrac{\G \g_J}{\G_2 \g_J} = H_1 \left( B_n; \Z[G] \right)$),
$$\dfrac{\G \g_K}{\G_2 \g_K} = H_1 \left( \EK; \Z[G] \right) \to H_1 \left( B_n; \Z[G] \right) = \dfrac{\G \g_J}{\G_2 \g_J}$$
is an isomorphism.\\

Consider the following commutative diagram of five-term exact sequences for $j < n$, with vertical maps induced by $\tau$:
$$
\xymatrix{
H_2\left( \HK \right) \ar[r]  & H_2\left( \dfrac{\HK}{\HKj} \right) \ar[r] \ar[d] & \dfrac{\HKj}{\HKjj} \ar[r] \ar[d] & 0 \\
H_2\left( \HJ \right) \ar[r] & H_2\left( \dfrac{\HJ}{\HJj} \right) \ar[r]  & \dfrac{\HJj}{\HJjj} \ar[r]  & 0 
}
$$
The image of $H_2\left( \HK \right)$ in $H_2\left( \dfrac{\HJ}{\HJj} \right)$ is precisely the image of $H_2\left( \HJ \right)$ by Proposition \ref{basics of \EK}.  We showed that $\dfrac{\HK}{\G_2 \gK} \cong \dfrac{\HJ}{\G_2 \gJ}$.  Suppose that $\dfrac{\HK}{\G_j \gK} \cong \dfrac{\HJ}{\G_j \gJ}$ for $j< n$.  Then $\dfrac{\HKj}{\HKjj} \cong \dfrac{\HJj}{\HJjj}$ by the commutativity of the above diagram.  It follows from the Five Lemma applied to 
$$
\xymatrix{
0 \ar[r] & \dfrac{\HKj}{\HKjj} \ar[r] \ar[d] & \dfrac{\HK}{\HKjj} \ar[r] \ar[d] & \dfrac{\HK}{\HKj} \ar[r] \ar[d] & 0 \\
0 \ar[r] & \dfrac{\HJj}{\HJjj} \ar[r] & \dfrac{\HJ}{\HJjj} \ar[r] & \dfrac{\HJ}{\HJj} \ar[r] & 0 
}
$$
that $\dfrac{\HK}{\HKjj} \cong \dfrac{\HJ}{\HJjj}$.  Applying the Five Lemma to
$$
\xymatrix{
0 \ar[r] & \dfrac{\HK}{\HKjj} \ar[r] \ar[d] & \dfrac{\PK}{\HKjj} \ar[r] \ar[d] & \dfrac{\PK}{\HKj} \ar[r] \ar[d] & 0 \\
0 \ar[r] & \dfrac{\HJ}{\HJjj} \ar[r] & \dfrac{\PJ}{\HJjj} \ar[r] & \dfrac{\PJ}{\HJj} \ar[r] & 0 
}
$$
we conclude that that $\dfrac{\PK}{\HKjj} \cong \dfrac{\PJ}{\HJjj}$ for all $j<n$.
\end{proof}
$\phantom{}$ \\

\subsection{indeterminacy of $\ka$-invariants}
\label{indeterminacy of ka-invariants}
We now investigate the dependence of $\ka_n(K, \tau)$ on $\beta$ and $\tau$.

\begin{prop}
\label{ka invariants and basing}
Suppose that $\tau_0: \EK \to B_n$ is an extendable map for the knot $K_{\beta_0}$.  If $\beta_1$ is another basing for $K$ then there is a canonical choice, up to based homotopy, of extendable map $\tau_1: \EK \to B_n$ for $K_{\beta_1}$ such that $[\ka_n(K_{\beta_0},\tau_0)]$ and $[\ka_n(K_{\beta_1},\tau_1)]$ differ in $[M, X_n]_0$ by an inner automorphism of $G$.  In particular, $\ka_n(K_{\beta_1},\tau_1)$ is freely homotopic to $\ka_n(K_{\beta_0},\tau_0)$. 
\end{prop}
\begin{proof} 
Let $K_i$ denote $K_{\beta_i}$ for $i=1, 2$, let $H_i = \pi_1(\EK \cup \beta_i, p)$, and let $\mu_i$ be the meridian of $K_i$.  Suppose that $\tau_0: E_{K_0} \to B_n$ is extendable and let $\tau_*$ be the homomorphism of fundamental groups induced by $\tau_0$.  If $a = [\beta_1 \ast \beta_0^{-1}]$ then conjugation by $\tau_*(a)$ induces an isomorphism $\sigma_{\tau_*(a)}: \pi_1\left( B_n \right) = \dfrac{\PJ}{\G_n \gJ} \to \dfrac{\PJ}{\G_n \gJ}$ such that $\sigma_{\tau_*(a)}(\mu_J) = \mu_J$ and the following diagram commutes:
$$
\xymatrix{
H_J \ar[d]^{\dagger_J} \ar[r]^{\cong}& H_J \ar[d]^{\dagger_J} \\
\dfrac{\PJ}{\G_n \gJ} \ar[r]^{\sigma_{\tau_*(a)}} & \dfrac{\PJ}{\G_n \gJ} 
}
$$
Since $B_n$ is an Eilenberg-MacLane space, $\sigma_{\tau_*(a)}$ induces a based homotopy equivalence $f: B_n \to B_n$ which is freely homotopic to the identity map.  Moreover, we may assume that $f$ restricts to a homeomorphism from $T_n$ to $T_n$ with $f(\mu_J) = \mu_J$ because $T_n$ and $\mu_J$ are also Eilenberg-MacLane spaces and the inclusions $\mu_J \hookrightarrow T_n \hookrightarrow B_n$ are cofibrations.  By an obstruction theory argument we may extend $f$ to a map $F: X_n \to X_n$ that is freely homotopic to the identity map.  The induced automorphism of $G$, $F_*$, is conjugation by $\gJ\left( \tau_*(a)\right)$.  Define $\tau_1 := f \circ \tau_0$, so $\tau_1: E_{K} \to B_n$ is an extendable map for $K_1$ and $\ka(K_1, \tau_1) = F \circ \ka(K_0, \tau_0)$.  Since $F$ is freely homotopic to the identity on $X_n$, $\ka_n(K_{1},\tau_1)$ is freely homotopic to $\ka_n(K_{0},\tau_0)$. 
\end{proof}

\begin{defn}
Define $A_n$ to be the group of automorphisms of the boundary condition
$ \left( \dagger_n, \mJ \right)$ over $G$,
$$\xymatrix{ A_n := \left\{ p: \left( \dagger_n, \mJ \right) \right. \ar[r]^(.48){\cong} & \left. \left( \dagger_n, \mJ \right) \Big| \gJ \circ p = \gJ \right\} .}$$
\end{defn}

\begin{prop}
\label{indeterminacy of kappa}
If $f,g : \EK \rightarrow B_n$ are extendable maps for a based knot $K_{\beta}$ then there is an element $p \in A_n$ inducing a based homotopy equivalence $P: X_n \to X_n$ such that $P_*: \pi_1(X_n) \to \pi_1(X_n)$ is the identity homomorphism and $P \circ \ka_n(K, f)$ is based homotopic to $\ka_n(K, g)$.  In particular, $P$ induces a bijection $P: [M, X_n]_0 \to [M, X_n]_0$ such that $P\left( [\ka_n(K, f)] \right) = [\ka_n(K, g)]$.
\end{prop}
\begin{proof}
By Proposition \ref{ext. induces isom.}, $f$ and $g$ induce isomorphisms $f_*, g_* : \dfrac{\PK}{\G_n \gK} \cong \dfrac{\PJ}{\G_n \gJ}$ over $G$ with $f_*(\mK) = g_*(\mK) = \mJ$ that restrict to isomorphisms between $\dagger_K\left( H_K \right)$ and $\dagger_n(H_J)$.  The automorphism $p := g_* \circ f_*^{-1} : \dfrac{\PJ}{\G_n \gJ} \rightarrow \dfrac{\PJ}{\G_n \gJ}$ is therefore an element of $A_n$.  By Proposition \ref{map of class space is WD}, $p$ induces a based homotopy equivalence $P: X_n \to X_n$, and hence a bijection $P: [M, X_n]_0 \to [M, X_n]_0$.  Since $P|_{B_n} \circ f$ is based homotopic to $g$, $P\left( [\ka_n(K, f)] \right) = [ P \circ \ka_n(K, f) ] = [\ka_n(K, g)]$.  By the construction of $P$ in the proof of Proposition \ref{map of class space is WD}, $P_*: \pi_1(X_n) \to \pi_1(X_n)$ is the identity homomorphism.
\end{proof}

\subsection{$n$-concordance}
\label{$n$-concordance sect.}
We define $n$-concordance and show that some $n^{th}$ $\ka$-invariant is defined and trivial if and only if $K$ is $n$-concordant to $J$.  As a consequence, if $K$ is concordant to $J$ then $K$ admits trivial $\ka$-invariants for all $n$.

\begin{defn}
\label{n-conc}
We say that a based knot $K_{\beta}$ is {\em $n$-concordant} to the based knot $J_{\alpha}$ if both of the following hold:

\begin{enumerate}
\item There is a cobordism $V \subset M \times I$ from $K$ to $J$.  Note that the inclusion of $i_V: \EV \rightarrow M\times I$ induces an epimorphism $\g_V: \PV = \pi_1(\EV) \rightarrow \M$, and we consider $\EV$ as a space equipped with this coefficient system.

\item When $V$ is based via $\alpha$ the composition 
$$\pi_1(\partial \EJ) \rightarrow \pi_1(\partial E_V) \rightarrow \PV \rightarrow \dfrac{\PV}{\HVn}$$
is an epimorphism onto the image of $\pi_1 (\partial E_V)$ in $\dfrac{\PV}{\HVn}$.
\end{enumerate}

We say that $K_{\beta}$ is {\em based $n$-concordant} to $J_{\alpha}$ if there is an $n$-concordance $V$ such that $\mu_{K_{\beta}} = \mu_{J_{\alpha}}$ in $ \dfrac{\PV}{\HVn}$. \\
\end{defn}

\begin{remark} 
\label{n-conc. remark}
$\phantom{}$

\begin{enumerate}
\item Unless stated otherwise, we regard $V$ as being based via $\alpha$ and we take $\mu_V = \mJ \times \{0\}$.  

\item Changing the basing of an $n$-concordance $V$ induces an inner automorphism of $\PV$.  In particular, there is an inner automorphism $\sigma: \PV \to \PV$ such that $\sigma
(\mJ) = \mK$.

\item For a {\em based} $n$-concordance we do not require that $\mK = \mJ$ in $\PV$; we only require that they be equal in $\dfrac{\PV}{\HVn}$. \\
\end{enumerate}
\end{remark}

\begin{prop}
\label{n-concordance has isom. n-quotient} 
If $\EV$ is an $n$-concordance between $\EK$ and $\EJ$ then the inclusion $\iota_L: E_L \to \EV$, for $L = K, J$, induces an isomorphism over $G$
$$ \dfrac{\pi_L}{\G_{j+1} \g_L} \cong \dfrac{\PV}{\HVjj} $$
for all $j\le n$. 
\end{prop}

\begin{proof}
We will prove the result for $K$ (the proof for $J$ is analogous) and so we regard $V$ as being based via the basing of $K$, since this makes sense when considering the inclusion $\iota = \iota_K$.  As in the proof of Proposition \ref{ext. induces isom.} we show that $\iota$ induces an isomorphism 
$$ \dfrac{\HK}{\G_2 \gK} \to \dfrac{\HV}{\G_2 \gV} $$ 
and then prove the general result by induction. \\

By construction, $M = \EK \cup_{\partial \EK} \overline{N(K)}$ and $M\times I = \EV \cup_{\partial \EV} \overline{N(V)}$ such that $M \hookrightarrow M \times I$ is a homology equivalence with any coefficients, $\partial \EK \hookrightarrow \partial \EV$, given by inclusion, induces a homeomorphism with one end of  $\partial \EV$, and $\overline{N(K)} \hookrightarrow \overline{N(V)}$ induces a homeomorphism with one end of $\overline{N(V)}$.  Consider the following diagram of long exact sequences,

$$
\xymatrix{
{} \ar[r] & H_2 \left( M; \Z[G] \right) \ar[r] \ar[d]^{\cong} & H_2 \left( M, \EK; \Z[G] \right) \ar[r] \ar[d] & H_1 \left( \EK; \Z[G] \right) \ar[r] \ar[d] & 0 &{}\save[]+<0cm,-1cm> *{\text{ $(\bigstar \bigstar)$}} \restore \\
{} \ar[r] & H_2 \left( M \times I; \Z[G] \right) \ar[r] & H_2 \left( M \times I, \EV; \Z[G] \right) \ar[r] & H_1 \left( \EV; \Z[G] \right) \ar[r] & 0 
}
$$
where $H_1\left(M; \Z[G] \right) = H_1\left(M \times I; \Z[G] \right) = 0$, being the homology of the universal covers of $M$ and $M \times I$, respectively.  By excision, 
$$H_2 \left( M, E_K; \Z[G] \right) \cong H_2 \left( \overline{N(K)}, \partial N(K); \Z[G] \right)$$
and
$$H_2 \left( M\times I, E_V; \Z[G] \right) \cong H_2 \left( \overline{N(V)}, \partial N(V); \Z[G] \right).$$
From the inclusion $\overline{N(K)} \hookrightarrow \overline{N(V)}$ we obtain
$$
\xymatrix{
0 \ar[r] & H_2 \left( \overline{N(K)}, \partial N(K); \Z[G] \right) \ar[r] \ar[d] & H_1 \left( \partial N(K); \Z[G] \right) \ar[r] \ar[d] & H_1 \left( \overline{N(K)}; \Z[G] \right) \ar[r] \ar[d] & 0 \\
0 \ar[r] & H_2 \left( \overline{N(V)}, \partial N(V); \Z[G] \right) \ar[r] & H_1 \left( \partial N(V); \Z[G] \right) \ar[r] & H_1 \left( \overline{N(V)}; \Z[G] \right) \ar[r] & 0
}
$$
where $H_2 \left( \overline{N(K)}; \Z[G] \right) = H_2 \left( \overline{N(V)}; \Z[G] \right) = 0$ since $\overline{N(K)}$ and $\overline{N(V)}$ each have the homotopy type of a 1-complexe.  If $V \times D^2 \cong \overline{N(V)}$ is a framing such that $\{ \ast \} \times S^1 \cong \mK$ for some $* \in K \subset V$ then, under this framing, the above commutative diagram becomes the following:
$$
\xymatrix{
0 \ar[r] & H_2 \left( \overline{N(K)}, \partial N(K); \Z[G] \right) \ar[r] \ar[d] & \txt{$H_1(K ; \Z[\M])$\\ $\oplus$\\ $H_1(\mK; \Z[P])$} \ar[r] \ar[d] & H_1 \left( K \times D^2; \Z[G] \right) \ar[r] \ar[d] & 0 \\
0 \ar[r] & H_2 \left( \overline{N(V)}, \partial N(V); \Z[G] \right) \ar[r] & \txt{$H_1(V ; \Z[\M])$\\ $\oplus$\\ $H_1(\mK; \Z[P])$} \ar[r] & H_1(V \times D^2; \Z[\M]) \ar[r] & 0
}
$$
where $P = \frac{\M}{C_{[K]}}$, the set of right cosets of the cyclic subgroup $C_{[K]} \leq G$ generated by $[K]$, and the coefficients in $H_1(\mK; \Z[P])$ are untwisted.  Since the second rightmost horizontal homomorphisms are induced by inclusion of $K \times S^1$ and $V \times S^1$ into the boundary of $K \times D^2$ and $V \times D^2$, respectively,
$$ H_1 (K; \Z[\M]) \to H_1(K \times D^2; \Z[\M]) $$
and 
$$ H_1 (V; \Z[\M]) \to H_1(V \times D^2; \Z[\M]) $$ 
are isomorphisms.  Moreover, $\mu_K$ bounds in $K \times D^2 \subset V \times D^2$, so these short exact sequences are split.  It follows that
$$H_2 \left( \overline{N(K)}, \partial N(K); \Z[G] \right) \cong H_1(\mK; \Z[P]) \cong H_2 \left( \overline{N(V)}, \partial N(V); \Z[G] \right).$$
Returning to the commutative diagram of long exact sequences $(\star \star)$, we see from the Five lemma that
$$\dfrac{\HK}{\G_2 \gK} = H_1 \left( \EK; \Z[G] \right) \cong H_1 \left( \EV; \Z[G] \right) = \dfrac{\HV}{\G_2 \gV}.$$\\

We want to proceed exactly as in the end of the proof of Proposition \ref{ext. induces isom.}.  For that, in the commutative diagram
$$
\xymatrix{
H_{2}(\HK) \ar[r] \ar[d]^{\iota_*} & H_{2}\left(\dfrac{\HK}{\HKj}\right) \ar[r] \ar[d]^{\iota_*}  & \dfrac{\HKj}{\HKjj} \ar[r] \ar[d]^{\iota_*}  & 0 \ar[d]\\
H_{2}(\HV) \ar[r] & H_{2}\left(\dfrac{\HV}{\HVj}\right) \ar[r] & \dfrac{\HVj}{\HVjj} \ar[r] & 0
}
$$
we need to see that $H_{2}(\HK)$ and $H_{2}(\HV)$ have the same image in $H_{2}\left(\dfrac{\HV}{\HVj}\right)$.  From the previous discussion we have a commutative diagram of exact sequences:
$$
\xymatrix{
0 \ar[r] & H_2 \left( \partial \EK; \Z[\M] \right) \ar[r] \ar[d] & H_2 \left( \EK; \Z[\M] \right) \ar[r]^{a_K} \ar[d] &H_2 \left( M; \Z[\M] \right) \ar[r] \ar[d]^{\cong} & H_1 \left( \mu_K; \Z[P] \right) \ar[d]^{\cong} \\
0 \ar[r] & H_2 \left( \partial \EV; \Z[\M] \right) \ar[r] & H_2 \left( \EV; \Z[\M] \right) \ar[r]^{a_V} &H_2 \left( M\times I; \Z[\M] \right) \ar[r] & H_1 \left( \mu_K; \Z[P] \right)
}
$$
If $M$ is irreducible or if $M \cong S^1 \times S^2$ and $[K]$ has infinite order in $\pi_1(S^1 \times S^2)$ then $a_K$ and $a_V$ are both the zero map (in the second case this is proved as in the proof of Proposition \ref{basics of \EK} (2)).  Hence, 
$$
\xymatrix{
H_{2}(\partial \EK; \Z[\M]) \ar@{->>}[r] \ar[d] & H_{2}(\EK;\Z[\M]) \ar@{->>}[r] \ar[d] & H_{2}(\HK) \ar[r] \ar[d] & H_{2}\left( \dfrac{\HK}{\HKj} \right) \ar[d] \\
H_{2}(\partial \EV; \Z[\M]) \ar@{->>}[r] & H_{2}(\EV;\Z[\M]) \ar@{->>}[r] & H_{2}(\HV) \ar[r] & H_{2}\left( \dfrac{\HV}{\HVj} \right)
}
$$ 
is commutative, with two epimorphisms in each row, as indicated.  It follows from the definition of $n$-concordance that $H_{2}(\partial \EK; \Z[\M])$ and $H_{2}(\partial \EV; \Z[\M])$ have the same image in $H_{2}\left(\dfrac{\HV}{\HVj}\right)$, as desired.  If $M \cong S^1 \times S^2$ and $K$ is null-homotopic then $a_K$ and $a_V$ are onto (see the proof of Proposition \ref{basics of \EK} (1)).  As in the proof of Proposition \ref{basics of \EK} (4), $H_2 \left( \partial \EK; \Z[\M] \right) \to H_{2}(\HK)$ and $H_2 \left( \partial \EV; \Z[\M] \right) \to H_{2}(\HV)$ are epimorphisms.\\

We now proceed exactly as in the end of the proof of Proposition \ref{ext. induces isom.}.  It was shown earlier that $\iota: \EK \hookrightarrow \EV$ induces an isomorphism $\dfrac{\HK}{\G_2 \gK} \cong \dfrac{\HV}{\G_2 \gV}$.  Suppose by induction that for some $j\leq n$, $\iota$ induces an isomorphism $\dfrac{\HK}{\HKj} \cong \dfrac{\HV}{\HVj}$.  From the above diagram it follows that $\dfrac{\HKj}{\HKjj} \cong \dfrac{\HVj}{\HVjj}$.  Applying the Five Lemma to 
$$
\xymatrix{
0 \ar[r] & \dfrac{\HKj}{\HKjj} \ar[r]\ar[d]_{\cong}^{\iota_*} & \dfrac{\HK}{\HKjj} \ar[r]\ar[d]^{\iota_*} & \dfrac{\HK}{\HKj} \ar[r]\ar[d]_{\cong}^{\iota_*} & 0\\
0 \ar[r] & \dfrac{\HVj}{\HVjj} \ar[r] & \dfrac{\HV}{\HVjj} \ar[r] & \dfrac{\HV}{\HVj} \ar[r] & 0\\
}
$$ 
we conclude that $\iota$ induces an isomorphism $\dfrac{\HK}{\HKjj} \rightarrow \dfrac{\HV}{\HVjj}$.  Applying the Five Lemma to
$$
\xymatrix{
0 \ar[r] & \dfrac{\HK}{\HKjj} \ar[r]\ar[d]_{\cong}^{\iota_*} & \dfrac{\PK}{\HKjj} \ar[r]\ar[d]^{\iota_*} & \M \ar[r]\ar[d]_{\cong}^{\iota_*} & 0\\
0 \ar[r] & \dfrac{\HV}{\HVjj} \ar[r] & \dfrac{\PV}{\HVjj} \ar[r] & \M \ar[r] & 0\\
}
$$ 
we see that $\iota$ induces an isomorphism $\dfrac{\PK}{\HKjj} \rightarrow \dfrac{\PV}{\HVjj}$ for all $j\leq n$.
\end{proof}

\begin{cor}
\label{n-conc is same as based}
If $K_{\beta}$ is $n$-concordant to $J_{\alpha}$ then there is some basing $\beta'$ such that $K_{\beta'}$ is based $n$-concordant to $J_{\alpha}$.
\end{cor}
\begin{proof}
Let $V$ be an $n$-concordance between $K_{\beta}$ and $J_{\alpha}$.  Then $\dfrac{ \pi_{K} }{\HKn} \cong \dfrac{ \pi_V }{\HVn}$ by Proposition \ref{n-concordance has isom. n-quotient}.  Since changing the basing of $V$ corresponds to an inner automorphism of $\PV$, there is an element $g \in \PK$ such that $g^{-1} \mK g \mapsto \mJ$ under $\PK \to \PV \to \dfrac{\PV}{\HVn}$.  Let $c$ be a simple closed curve representing $g^{-1}$ and define $\beta' = c \ast \beta$.  Then $V$ is a based $n$-concordance between $K_{\beta'}$ and $J_{\alpha}$.
\end{proof}

For the next proposition recall from Example \ref{can. extension} that the canonical epimorphism $q_n: \PJ \to \dfrac{\PJ}{\HJn}$ induces an extendable map $\tau_n: \EJ \to B_n$ up to based homotopy.
\begin{prop}
\label{n-concordance gives extendable}
Let $V$ be an $n$-concordance.  For all $j\leq n$ there exists an extendable map $f: \EV \rightarrow B_j$.  Moreover, the restriction of $f$ to $\EJ \times \{0\}$ is based homotopic to $\tau_j$.
\end{prop}
\begin{proof}
By Proposition \ref{n-concordance has isom. n-quotient} the inclusion $\iota = \iota_J : \EJ \hookrightarrow \EV$ induces an isomorphism $\iota_*: \dfrac{\PJ}{\HJj} \cong \dfrac{\PV}{\HVj}$ for $j\leq n+1$.  We therefore get a homomorphism $f_*: \PV \rightarrow \dfrac{\PJ}{\HJj}$ as in the following commutative diagram:
$$
\xymatrix{
\PJ \ar[r]^{q_j} \ar[d]_{\iota_*} & \dfrac{\PJ}{\HJj} \ar[d]^{\iota_*}_{\cong} \\ 
\PV \ar[r] \ar@{..>}[ur]^{f_*} & \dfrac{\PV}{\HVj}  .
}
$$ 
As $B_j$ is an Eilenberg-MacLane space, $f_*$ induces a map $f: \EV \rightarrow B_j$ up to based homotopy.  Since $V$ is an $n$-concordance,  $\mu_V = \mJ \times \{0\}$ and
$$
\xymatrix{ 
\pi_1(\partial \EJ) \ar[r] & \pi_1(\partial \EV) \ar[r] & \PV \ar[r] & \dfrac{\PV}{\HVn} \ar[r] & \dfrac{\PV}{\HVj} 
}
$$ \\ 
is an epimorphism onto the image of $\pi_1(\partial \EV)$ in $\dfrac{\PV}{\HVj}$ for $j\leq n$.  In particular, 
$$f_*\left( \pi_1(\partial \EV \cup \alpha, p\times 0) \right) = \dagger_j(H_J).$$ 
The inclusion $\partial \EV \hookrightarrow \EV$ is a cofibration, so we may take $f: \EV \to B_j$ such that $f(\partial \EV) \subseteq T_j$ and $f(\mu_V) = \mJ$.  Finally, $f_* \circ \iota_* = q_j$, so $f$ restricted to $\EJ \times \{0\}$ is based homotopic to $\tau_j$.  It follows that $f$ is the desired map.  
\end{proof}

\begin{thrm}
\label{n-conc iff nth invt}
A knot $K$ is $n$-concordant to $J$ if and only if some $j^{th}$ $\ka$-invariant $\ka_j(K, \tau)$ is defined and trivial for all $j\leq n$.
\end{thrm}
\begin{proof}
Suppose that $K$ is $n$-concordant to $J$ via $V$.  We may assume by Corollary \ref{n-conc is same as based} and Proposition \ref{ka invariants and basing} that $V$ is a based $n$-concordance.  By Proposition \ref{n-concordance gives extendable} there exists for each $j\leq n$ an extendable map $f_j: \EV \rightarrow B_j$ that when restricted to $\EJ \times \{0\}$ is based homotopic to $\tau_j$.  Restricting $f_j$ to $\EK \times \{1\}$ defines an extendable map $\EK \rightarrow B_j$, so a $j^{th}$ $\ka$-invariant $\ka_j(K, {f_j})$ of $K$ exists.  Extending $f_j$ to $F_j : M\times I \rightarrow X_n$ defines a homotopy from $\ka_j(K, {f_j})$ to $\ka_j(J, \tau_j)$, so $[\ka_j(K, {f_j})]$ is trivial.\\

Suppose now that an $n^{th}$ $\ka$-invariant $\ka_n(K, \tau)$ of $K$ is defined and trivial.  Then there is a homotopy 
$$\xymatrix{ H:M\times I \ar[r] & X_n }$$ 
with 
\begin{equation*}
H(-, s) = \left\{
\begin{array}{ccl}
\ka_n(J, \tau_n)(-) & \text{if} & s = 0 \\
\ka_n(K, \tau)(-) & \text{if} & s = 1 \\
\end{array} \right.
.
\end{equation*} 
Choose a framing $t: S^1 \times D^2 \cong ST$ with $\{ \ast \} \times \partial D^2 \cong \mJ$ for some $\ast \in S^1$ and let $S = t(S^1 \times \{0\})$.  We may assume that $H$ is transverse to $ST$ and contains $S \subset ST$ in its image.  Pulling back $S$ yields a surface 
$$\xymatrix{ V = V_0 \coprod \cdots \coprod V_m }$$ 
in $M\times I$, transverse to the boundary, such that $V\cap \left(M\times \{0\} \right) = J\times \{0\}$ and $V\cap \left(M\times \{1\} \right) = K\times \{1\}$, and each $V_i$ is connected.  If $V$ is connected (and based via the basing of $J$) then $H|_{\EV}: \EV \rightarrow B_n$ is an extendable map and $V$ is an $n$-concordance. \\

Suppose that $V$ is not connected.  We now show how to connect the disjoint pieces of $V$ and redefine $H$ such that the new map is extendable, thus obtaining an $n$-concordance from $K$ to $J$.  Let $q \in S$ and for $i \in \{0, 1\}$ let $q_i \in V_i$ be points such that $H(q_i) = q$.  Let $B(q) \cong I \times D^2$ be a closed regular neighborhood of a subinterval of $S$ containing $q$ with $I\times \{0\} \subset I\times D^2$ a zero-section.  Let $B(q_i) \cong D^2 \times D^2$ be a regular neighborhood of a closed disk in $V_i$ containing $q_i$ such that $H(B(q_i)) = B(q)$ and such that $D^2 \times \{0\} \subset D^2 \times D^2$ is the zero-section.  Attach a one-handle to $(M \times I) \times I$ along $(M\times I)\times \{1\}$ to obtain a new manifold 
$$\xymatrix{ N = \left(M\times I \right) \times I \cup D^2 \times D^2 \times I. }$$ 
We can arrange that $D^2 \times D^2\times \{i\}$ is attached to $B(q_i)\times \{1\}$ such that
\begin{enumerate}
\item $D^2 \times \{0\} \times \{i\}$ is identified with the zero-section of $B(q_i) \times \{1\}$,
\item The resulting 3-dimensional submanifold $W_1 = V \times I \cup D^2 \times \{0\} \times I$ is orientable, and
\item $N$ is orientable.
\end{enumerate}
We give $N$ the orientation of $(M\times I)\times I$. \\

Since 
$$
\xymatrix{
\EJ \ar@{^(->}[dr] \ar[rr] & & B_n \\
& E_V \ar[ur]
}
$$ 
commutes, $H_*: \pi_1(E_V \times I, p\times 0 \times 1) \rightarrow \dfrac{\PJ}{\HJn}$ is surjective.  Define meridians $\mu_0$ and $\mu_1$ for $V_0$ and $V_1$, respectively, as follows:  Choose a path $\alpha_i$, for $i \in \{0, 1 \}$, from $p\times 0 \times 1$ to a point $v_i \times e_i \in D^2 \times S^1 \subset B(q_i)$ and let $\mu_i = \{v_i\} \times S^1 \cup \alpha_i$.  Since $H_*$ is surjective, we can choose $\alpha_i$ such that $H_*(\mu_i) = \mJ$.  Recall that 
$$ W_1 = V \times I \cup D^2 \times \{0\} \times I ,$$
and let $N(W_1)$ be an open regular neighborhood of $W_1$.  Then
$$\pi_1\left( N - N(W_1) \right) = \dfrac{ \pi_1 \left( M\times I - N(V \times I) \right) \ast \Z }{ \langle \mu_0 \sim a^{-1} \mu_1 a \rangle },$$
where $a$ generates $\Z$ and if $y:I \to N$ runs along the one-handle from $v_0 \times e_0$ to $v_1 \times e_1$, without intersecting $(D^2 \times \{0\}) \times I$, then $a = [\alpha_1 \ast y^{-1} \ast \alpha_0^{-1}]$.  Define an epimorphism 
$$h: \pi_1\left( N - N(W_1) \right) \to \dfrac{\PJ}{\HJn}$$
over $\M$ by defining $h(a)$ to be trivial (that is, $h(a)$ is the identity element) and defining $h$ to be $(H|_{M\times I - N(V \times I)})_*$ on $\pi_1 \left( M\times I - N(V \times I) \right)$.  This is well-defined since we chose $\mu_0$ and $\mu_1$ such that $H_*(\mu_0) = H_*(\mu_1) = \mJ$.  Since $H\left( \partial N(V \times I) \right) \subseteq T_n$ and $h(\mu_0) = h(\mu_1) = \mJ$, it follows that $h \left( \pi_1\left( \partial N(W_1) \right) \right) = \dagger_n(H_J)$, where we base $W_1$ via the basing of $(J \times \{0\} ) \times \{1 \}$.  Since $B_n$ is an Eilenberg-MacLane space, $h$ extends $H$ to a based map $\overline{H}: N - N(W_1) \to B_n$ such that $\overline{H} \left( \partial N(W_1) \right) \subseteq T_n$.  In particular, $\overline{H}$ restricted to $N - N(W_1)$ is an extendable map to $B_n$, after choosing suitable basings for the components of $W_1$ (which can be done since, as we observed earlier, $H_*$ is surjective). \\

Our goal now is to add a two-handle to cancel the previously attached one-handle.  
Let 
$$c_{0,1}: I \rightarrow E_V \times I \subset (M\times I) \times I$$ 
be an embedded path from $\partial B(q_0) \times \{1\}$ to $\partial B(q_1) \times \{1\}$ such that $\overline{H}\left( c_{0,1}(0) \right)$ and $\overline{H}\left( c_{0,1}(1) \right)$ lie in the image of the basing path $\alpha$ in $B_n$.  Recall that $H_*: \pi_1(E_V \times I, p\times 0 \times 1) \rightarrow \dfrac{\PJ}{\HJn}$ is surjective.  We may therefore assume that $\overline{H}(c_{0,1}) \cup \alpha$, the path in $B_n$ given by following $\alpha$ until it reaches $\overline{H}\left( c_{0,1}(0) \right)$, following $\overline{H}(c_{0,1})$ until it reaches $\overline{H}\left( c_{0,1}(1) \right)$, and then going backward down $\alpha$, is homotopically trivial in $B_n$.  Let 
$$c:I \rightarrow E_{W_1} \subset N$$
be a curve that travels from $c_{0,1}(0) \in \partial B(q_0) \times \{1\}$ to $c_{0,1}(1) \in \partial B(q_1) \times \{1\}$ through the one-handle and geometrically intersects the belt of the one-handle once.  Since $h: \pi_1\left( N - N(W_1) \right) \to \dfrac{\PJ}{\HJn}$ is surjective and $h(a)$ is trivial, we may assume that $\overline{H}(c) \cup \alpha$ is homotopically trivial in $B_n$.  Then $c_{0, 1} \ast c^{-1}$ is a closed embedded curve in $N$ and $\overline{H}(c_{0, 1} \ast c^{-1}) \cup \alpha$ is homotopically trivial in $B_n$.  We may isotope $c_{0, 1} \ast c^{-1}$ to lie embedded on the boundary of $N$, still geometrically intersecting the belt of the previously attached one-handle exactly once, and then attach a two-handle to cancel this one-handle.  Since $c_{0, 1} \ast c^{-1}$ is null-homotopic in $B_n$, $\overline{H}$ extends over this two-handle.  The result is a manifold diffeomorphic to $(M\times I)\times I$ together with a submanifold $W_1$ with one less connected component than $V \times I$.  Moreover, $\overline{H}$ restricts to an extendable map on the exterior of $W_1$, after an appropriate choice of basings for its components. \\

After connecting the components of $V$ we will have a cobordism between two knots, but we still need to see that those two knots are $J$ and $K$.  However, we can choose the diffeomorphism between $(M\times I)\times I$ and the manifold we obtained by adding the one- and two-handles to be the identity on $(M\times I) \times \{0\}$.  Hence, the knots $W_1 \cap (M\times \{0\}) \times \{1\}$ and $W_1 \cap (M\times \{1\}) \times \{1\}$ are isotopic to $J$ and $K$, respectively. \\


Let $W_m \subset (M\times I) \times I$ be the manifold obtained by repeating this procedure $m$ times to connect the components of $V = V_0 \coprod \cdots \coprod V_m$.  The connected surface $W_m \cap (M\times I) \times \{1\}$, based via the basing of $(J\times \{0\} ) \times \{1\}$ in $(M\times \{0\} ) \times \{1\}$, is our $n$-concordance.
\end{proof}

\begin{cor}
\label{ka is conc invt}
If $K$ is concordant to $J$ then for each $n$ there is an extendable map $\tau$ such that $\ka_n(K, \tau)$ is defined and trivial. 
\end{cor}
\begin{proof}
A concordance is an $n$-concordance for all $n$.
\end{proof}

\subsection{Characteristic knots}
\label{Characteristic knots sect.}
We now restrict our attention to knots in manifolds of the second type (i.e. closed, orientable, aspherical 3-manifolds) and investigate the relationship between $\ka$-invariants and $J$-characteristicness.  \\

Recall that a knot $K$ in $M$ is {\em $J$-characteristic} if there is a continuous degree-one map $\alpha: M \to M$ such that $\alpha(K) = J$ and $\alpha(M-K) \subseteq M-J$.

\begin{thrm}
\label{dense J}
Let $K$ be a knot in an orientable, aspherical 3-manifold.  If $K$ is $J$-characteristic via a map $\alpha:M \to M$ that indices the identity on $G$ then $K$ has trivial $\ka$-invariants.
\end{thrm}
\begin{proof}
Up to based homotopy, $\alpha$ induces a map $\tau: \EK \to \EJ$ that restricts to a homeomorphism $\partial \EK \cong \partial \EJ$ with $\tau(\mu_K) = \mu_J$.  Hence, the composition $\tau_n \circ \tau: \EK \to B_n$ is extendable for each $n$ such that the following diagram commutes:
$$
\xymatrix{
M \ar[dr]_{\ka(K, \tau_n \circ \tau)} \ar[rr]^{\alpha} & & M \ar[dl]^{\ka(J, \tau_n)} \\
& X_n
}
$$
Since $\alpha$ indices the identity on $G$ and $M$ is an Eilenberg-MacLane space, $\alpha$ is based homotopic to the identity map on $M$.  In particular, $\ka(K, \tau_n \circ \tau)$ is homotopic to $\ka(J, \tau_n )$.
\end{proof}

\begin{remark}
Let $K$ be a knot in $M$, let $\eta \subset \EK$ be an embedded curve that bounds an embedded disk in $M$, and let $N(\eta)$ be a regular neighborhood of $\eta$ in the interior of $\EK$.  Let $L$ be a knot in $S^3$.  There is a canonical choice of longitude for both $\eta$ and $L$, which we denote $\lambda_{\eta}$ and $\lambda_L$, respectively; let $\mu_{\eta}$ and $\mu_L$ be their respective meridians.  Then 
$$\left( M - N(\eta) \right) \cup -\left( S^3 - N(L) \right), $$
where $\mu_L \sim \lambda_{\eta}^{-1}$ and $\lambda_L \sim \mu_{\eta}$, is diffeomorphic to $M$ via a diffeomorphism that is the identity outside a regular neighborhood of the disk bounded by $\eta$.  In particular, this diffeomorphism induces the identity homomorphism on $G$.  
\end{remark}

\begin{defn}
Define $K(\eta, L)$ to be the image of $K$ in $\left( M - N(\eta) \right) \cup -\left( S^3 - N(L) \right) \cong M$.  That is, the image of $K$ under the inclusion
$$K \hookrightarrow M-N(\eta) \hookrightarrow \left( M - N(\eta) \right) \cup -\left( S^3 - N(L) \right).$$ 
\end{defn}

For a special case of the following lemma see Proposition 2.2 of \cite{L2-paper}.  The proof is the same.
\begin{lem}
\label{EMS}
If $K$ is a homotopically essential knot in a closed, orientable, aspherical 3-manifold then $\EK$ is an Eilenberg-MacLane space.
\end{lem}

\begin{thrm}
\label{char has triv}
Let $J$ be a homotopically essential knot in a closed, orientable, aspherical 3-manifold.  Let $\eta$ be a curve in $\EJ$ that bounds an embedded disk in $M$ and let $L$ be a knot in $S^3$.  Then $J(\eta, L)$ has trivial $\ka$-invariants.
\end{thrm}
\begin{proof}
Since 
$$E_{J(\eta, L)} = \left( \EJ - N(\eta) \right) \cup -\left( S^3 - N(L) \right), $$ 
where $\mu_{L} \sim \lambda_{\eta}^{-1}$ and $\lambda_{L} \sim \mu_{\eta}$, it follows from the Seifert-van Kampen theorem that 
$$\pi_{J(\eta, L)} := \pi_1(E_{J(\eta, L)}) =  \pi_1(\EJ - N(\eta)) \ast_{\Z^2} \pi_1(S^3 - N(L)),$$ 
as in the following pushout diagram,
$$
\xymatrix{
\Z^2 \ar[r] \ar[d] & \pi_1(S^3 - N(L)) \ar[d] \\
\pi_1(\EJ - N(\eta)) \ar[r] & \pi_{J(\eta, L)}
}
$$
where $\Z^2 = \pi_1(\partial N(\eta))$.  Also, the image of $\mJ$ in $\pi_{J(\eta, L)}$ is a meridian, $\mu_{J(\eta, L)}$.  Let $f:  \pi_1(\EJ - N(\eta)) \to \PJ$ be the epimorphism induced by inclusion $\EJ - N(\eta) \hookrightarrow \EJ$ and let $g: \pi_1(S^3 - N(L)) \to \Z \to \PJ$ be abelianization followed by the homomorphism that sends $\mu_{L}$ to $\lambda_{\eta}^{-1}$.  Then the pushout property of $\pi_{J(\eta, L)}$ induces a homomorphism $h: \pi_{J(\eta, L)} \to \PJ$ over $\M$ as in the following diagram:
$$
\xymatrix{
\Z^2 \ar[r] \ar[d] & \pi_1(S^3 - N(L)) \ar[d] \ar@/^/[ddr]^g \\
\pi_1(\EJ - N(\eta)) \ar[r] \ar@/_/[drr]_f & \pi_{J(\eta, L)} \ar@{..>}[dr]|h \\
& & \PJ 
}
$$
Since $f(\mJ) = \mJ$, the commutativity of the diagram implies that $h(\mu_{J(\eta, L)}) = \mJ$.  Moreover, the inclusion $\partial \EJ \hookrightarrow \EJ - N(\eta)$ induces a commutative diagram:
$$
\xymatrix{
{} & {} & \pi_{J(\eta, L)} \ar[dd]^h \\
H_J \ar[r] & \pi_1 \left( \EJ - N(\eta) \right) \txt{$\phantom{KKKK}$} \ar[ur] \ar[dr] & {}\\
{} & {} & \pi_J
}
$$
Because $\EJ$ is an Eilenberg-MacLane space by Lemma \ref{EMS} and $\mu_J \hookrightarrow \partial \EJ \hookrightarrow \EJ$ are cofibrations, $h$ induces a map $\tau: E_{J(\eta, L)} \to \EJ$ that extends to $M \to M$ and induces the identity on $G$.  As in the proof of Theorem \ref{dense J}, we see that $J(\eta, L)$ has trivial $\ka$-invariants.
\end{proof}

\bibliography{BigBibliography}
\bibliographystyle{plain}

\end{document}